\newcommand{\textbfME}{}
\newcommand{\beq}{\begin{equation}}
\newcommand{\eeq}{\end{equation}}
\newcommand{\beqa}{\begin{eqnarray}}
\newcommand{\eeqa}{\end{eqnarray}}
\newcommand{\bal}{\begin{align}}
\newcommand{\eal}{\end{align}}
\newcommand{\bfl}{\begin{flalign}}
\newcommand{\efl}{\end{flalign}}
\newcommand{\bbm}{\begin{bmatrix}}
\newcommand{\ebm}{\end{bmatrix}}

\newcommand{\pdd}[2]{\frac{\partial #1}{\partial #2}}

\newcommand{\sgn}{\operatorname{sgn}}

\newcommand{\boldsymbolME}{}
\newcommand{\bm}[1]{\boldsymbolME{#1}}
\newcommand{\vt}[1]{\textbfME{#1}}
\RequirePackage[fleqn]{amsmath}

\documentclass{article}       

\usepackage{amsmath}
\usepackage{bm}
\usepackage[normalem]{ulem} 
\usepackage{color}
\definecolor{darkgreen}{rgb}{0.0,0.6,0.0}

\usepackage{hyperref}
\usepackage{color}
\usepackage{amssymb}
\usepackage{comment}
\usepackage{graphicx}
\usepackage{fmtcount}
\usepackage{cite}
\usepackage[lofdepth,lotdepth]{subfig}
\usepackage{tikz}
\usetikzlibrary{decorations.pathmorphing,patterns}
\usepackage{pgfplots}
\pgfplotsset{compat=1.3}
\begin{document}

\title{Bifurcation analysis of a smoothed model of a forced impacting beam and comparison with an experiment
}

\author{M.\ Elmeg{\aa}rd\thanks{
              Department of Applied Mathematics and Computer Science, Technical University of Denmark, Matematiktorvet, Building 303B, 2800 Kgs. Lyngby, Denmark, melm@dtu.dk.
              }   \and
        B.\ Krauskopf\thanks{Department of Mathematics, University of Auckland, Private Bag 92019, Auckland 1142, New Zealand, b.krauskopf@auckland.ac.nz.} \and
        H.M.\ Osinga\thanks{Department of Mathematics, University of Auckland, Private Bag 92019, Auckland 1142, New Zealand, h.m.osinga@auckland.ac.nz.} \and
         J.\ Starke\thanks{Department of Applied Mathematics and Computer Science, Technical University of Denmark, Matematiktorvet, Building 303B, 2800 Kgs. Lyngby, Denmark, jsta@dtu.dk} \and
         J.J.\ Thomsen\thanks{Department of Mechanical Engineering, Technical University of Denmark, Nils Koppels Alle 404, DK-2800 Kongens Lyngby, Denmark, jjt@mek.dtu.dk} 
}

\maketitle

\vspace{-0.5cm}
\begin{abstract}
A piecewise-linear model with a single degree of freedom is derived from first principles for a driven vertical cantilever beam with a localized mass and symmetric stops. The resulting piecewise-linear dynamical system is smoothed by a switching function (nonlinear homotopy). For the chosen smoothing function it is shown that the smoothing can induce bifurcations in certain parameter regimes. These induced bifurcations disappear when the transition of the switching is sufficiently and increasingly localized as the impact becomes harder. The bifurcation structure of the impact oscillator response is investigated via the one- and two-parameter continuation of periodic orbits in the driving frequency and/or forcing amplitude. The results are in good agreement with experimental measurements.\\
\vspace{-0.5cm}
\end{abstract}

\newpage

\section{Introduction}

In mechanical engineering applications piecewise-smooth dynamical systems are often encountered through structural interactions. This is especially the case in the field of machinery dynamics, e.g.,\ with vibro-impact, and with friction systems and processes: \textit{Vibro-impact systems} and processes involve repeatedly colliding elements and vibrations with abrupt changes in velocity and forces. Applications include devices to crush, grind, forge, rivet, drill, punch, tamp, print, tighten, pile, cut, and surface treat a variety of materials and objects, at frequencies ranging from sub-Hertz to ultrasonic \cite{Babitsky1998}, \cite{Burton1968}, \cite{Kobrinskii1969}. In other contexts vibro-impact occurs as an unwanted side effect, often producing noise and wear, such as with devices operating with stops and clearances, e.g., gear wheels, rotors, rattling heat exchanger tubes, and guides for roller chain drives. Here, discontinuities are typically due to abrupt changes in the restoring forces or boundary conditions, i.e.\ on a time scale much smaller than those associated with the free oscillation frequencies of a system. \textit{Friction systems} and processes are an inherent or purposefully implemented part of many technical devices, including automotive brakes and clutches and bowed musical instruments, some are unwanted in other apparatus, such as in chatter in metal cutting, creaking doors, and squealing tramways and disc brakes \cite{RA1992a}, \cite{RA1992b}. Here, discontinuities appear in the dissipative terms of the equations of motion, typically as abrupt changes in dissipative forces with the sign changes of the relative interface sliding velocities.\par 
The strong nonlinearity associated with such discontinuous systems and processes preclude exact analytical solution in all but a few simple cases and, in general, bifurcation analysis is far from straightforward \cite{Leine2004}. Analytical approaches are typically approximate; they include ÒstitchingÓ \cite{Kobrinskii1969}, i.e.\ integrating motions between impacts or other discontinuities, and using kinematic impact conditions to switch solution intervals; averaging; harmonic linearization \cite{Babitsky1998}; and various kinds of discontinuous transformations of variables, e.g.,\ see\cite{Thomsen2008}. Details about some of the mathematical methods for piecewise-smooth systems can be found, e.g.,\ in\cite{bernardo2007piecewise}. However, in most cases one needs to resort to or complement the analysis with numerical solutions of the equations of motion \cite{Acary2008}.\par
One of the main modeling challenges emerges from the fact that the continuum physical laws of impacts, friction etc.\ are not well understood. Hence, while the equations of motion for many smooth mechanical systems can be derived directly from the theory of elasticity, the modeling of mechanical systems with impacts and/or friction involves identifying suitable reaction laws. This results in a range of modeling decisions that are often problematic to substantiate theoretically. This is both in terms of the mathematical structure of the reaction laws, as well as the extra parameters that need to be estimated. Under these conditions it is very advantageous to be able to make comparisons with experimental measurements to test the given model.\par
In this paper we construct a mathematical model for a vibro-impacting cantilever of which experiments were performed and reported in \cite{Bureau2013}; see Figure~\ref{fig:system_sketch} for a sketch of the experimental setup. In \cite{Bureau2013}, the authors perform control-based bifurcation analysis of the experiment, i.e.,\ real-time bifurcation analysis of a mechanical vibro-impacting system. The method of control-based continuation is a method that makes it possible to perform bifurcation analysis of systems for which no mathematical model is given; the method was recently developed in \cite{Sieber2008} and demonstrated in a rotating pendulum experiment in \cite{Sieber2008a}. Since then, control-based continuation has been used to investigate the bifurcation structures of mechanical experiments, such as the current case of a cantilever beam with symmetric mechanical stops \cite{Bureau2013} and a nonlinear energy harvester \cite{Barton2011}, \cite{Barton2012}. The model derived here deals with the smooth modeling of small-amplitude vibrations of a driven cantilever beam under the influence of symmetric stops, and it concludes with a comparison with experimental measurements. In Section \ref{sec:math_model} we derive the single-degree-of-freedom model of a cantilever beam that vibro-impacts on symmetric stops. The structural mechanics of the model is based on the Euler-Lagrange beam equation and the reduction is inspired by the Galerkin step in \cite{MOON1983} and \cite{SHAW1985}. In Section \ref{sec:smoothing} a smooth version of the dynamical system is constructed by a nonlinear homotopy. In Section \ref{sec:num_bif} numerical bifurcation analysis of the smoothing procedure is performed in order to test that the bifurcation results are robust with respect to the smoothing homotopy parameter. In Section \ref{sec:exp_comp} the model results are compared with experimental measurements.

\section{Mathematical model}\label{sec:math_model}
\begin{figure}[h]
\centering
\begin{picture}(5,10)
\put(0,0){\includegraphics[width=0.4\textwidth]{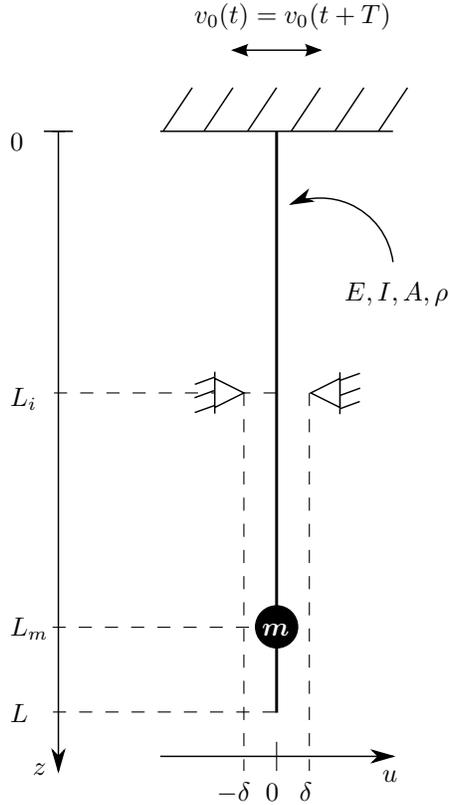}}
\put(-0.45,8.3){$0$}
\put(-0.45,4.9){$L_i$}
\put(-0.45,1.85){$L_m$}
\put(-0.45,0.7){$L$}
\put(-0.15,0){$z$}
\put(2.95,-0.35){$0$}
\put(2.3,-0.35){$-\delta$}
\put(3.4,-0.35){$\delta$}
\put(4.5,-0.1){$u$}
\put(2.9,1.85){\textcolor{white}{{$\boldsymbol{m}$}}}
\put(2.,10.){$v_0(t)=v_0(t+T)$}
\put(4,6.3){$E,I,A,\rho$}
\end{picture}
\caption{Experimental setup of the vertical cantilever beam of length $L$, with a lumped mass point $m$, symmetric mechanical stops at a transverse distance $\delta$ and a periodic excitation $v_0(t)$ of the clamped boundary. In the experiment; $E=2\cdot 10^{11}\, \rm Pa$, $I=2.08\cdot10^{-12}\, {\rm m}^4$, $A=2.5\cdot10^{-5}\, {\rm m}^2$, $\rho = 8\cdot 10^3\, \frac{\rm kg}{{\rm m}^3}$, $m = 0.2116\, {\rm kg}$, $\delta \approx 10^{-3}\, \rm m$, $L_m=0.1275\, \rm m$, $L_i=0.071\, \rm m$ and $L_m=0.1275\, \rm m$ and $L=0.161\, \rm m$.}
\label{fig:system_sketch}
\end{figure}

A mathematical model for the mechanical vibrations of a cantilever beam with symmetric mechanical stops, as sketched in Figure \ref{fig:system_sketch}, is derived. In Section~\ref{subsec:PDE} the governing PDE and the main assumptions for its validity are given. In Section \ref{subsec:PDE_to_ODE} the structure of the single-degree-of-freedom  model is derived by Galerkin's method, but without an explicit choice for the approximating subspace, i.e.,\ a mode shape is not chosen explicitly. In Section \ref{subsec:Parameter_estimation} the model parameters are estimated and compared with values fitted to experimental data.

\addtolength{\abovecaptionskip}{3mm}  
 
\subsection{Transverse vibration of a cantilever with a lumped mass} \label{subsec:PDE}

Formally the PDE is derived by Euler-Bernoulli beam theory, for small-amplitude transverse vibrations of a slender cantilever. Furthermore, the cantilever is harmonically driven and has a lumped mass $m$ attached. Following, e.g.,\ \cite{Landau1986}, \cite{JuelThomsen2003}, the complete system can be derived as
\begin{subequations}
\label{eq:PDE_all}
\begin{align}
&\Big(m\delta(z-L_m)+\rho A_0\Big)\ddot{u} +EI {u}''''+\mu \dot{u}''''=0,\label{eq:PDE}\\
&u'(0,t)=u''(L,t)=u'''(L,t)=0 \quad \text{ and } \quad u(0,t)=v_0(t).\label{eq:PDE_BC}
\end{align}
\end{subequations}
Here, $u=u(z,t)$ is the transverse displacement in inertial coordinates, $z$ is the axial coordinate; the overdot denotes derivation w.r.t.\ time $t$, and the prime denotes derivation w.r.t.\ space $z$. Furthermore, $\delta(\cdot)$ is the Dirac delta distribution; $v_0(t)=A\cos\omega t$ describes the harmonic excitation of the clamped base; $L$ is the beam length; $\rho$ is the beam density; $A_0$ is the cross-sectional area of the beam; $\mu$ is the coefficient of stiffness proportional damping; $E$ is the modulus of elasticity; $I$ is the cross section area moment of inertia; and $L_m$ is the axial position of the attached mass.
In a co-moving frame, using $ u=\eta+v_0$, the governing equations transform to 
\begin{subequations}
\label{eq:PDE_comoving_all}
\begin{align}
&\Big(m\delta(z-L_m)+\rho A_0\Big)(\ddot{\eta}+ \ddot{v}_0) + EI{\eta}''''+\mu\dot{\eta}'''' =0, \label{eq:PDE_comoving}\\
&\eta(0,t)=\eta'(0,t)=\eta''(L,t)=\eta'''(L,t)=0.\label{eq:PDE_comoving_BC}
\end{align}
\end{subequations}
The equations of motion in the co-moving frame \eqref{eq:PDE_comoving_all} are equivalent to a cantilever beam with a lumped mass under the influence of an axially uniformly distributed harmonic forcing, and a time-harmonic point force proportional to the lumped mass at the position $z=L_m$.

\subsection{Reduced ODE model}\label{subsec:PDE_to_ODE}

The configuration space of the cantilever in the PDE description is infinite dimensional. In order to construct a finite-dimensional description an approximating subspace is chosen for the spatial coordinate. Let $\eta(z,t) = \sum\limits_{i=1}^{N} a_i(t)\phi_i(z)$, where all $\phi_i(z)$ satisfy the essential boundary conditions. Galerkin's method is then used to obtain a finite-dimensional ODE 
\begin{align}
&\int_0^L \Big[\Big(m\delta(z-L_m)+\rho A_0\Big)(\ddot{\eta}+ \ddot{v}_0) + EI {\eta}''''+\mu\dot{\eta}''''\Big]\phi_j(z) \text{d}z=0. \label{eq:trans_beam_co-moving}
\end{align}
After applying essential and natural boundary conditions via integration by parts, Equation \eqref{eq:trans_beam_co-moving} can be written as
\begin{align}
\label{eq:secondorderlinear}
& \bold{M}\ddot{\bold{a}}+\bold{C}\dot{\bold{a}}+\bold{K}\bold{a}+\ddot{v}_0\bold{d}=0,
\end{align}
where
\begin{subequations}
\label{eq:parameters}
\begin{align}
M_{ij} &= m \phi_i(L_m)\phi_j(L_m)+\rho A_0\big<\phi_i,\phi_j\big>, \label{eq:mass_coefficient} \\
C_{ij} &= \mu \big<\phi_i'',\phi_j''\big>, \label{eq:damping_coefficient} \\
K_{ij} &= EI \big<\phi_i'',\phi_j''\big>, \label{eq:stiffness_coefficient}\\
d_{i} &= m \phi_i(L_m)+\rho A_0\big<1,\phi_i\big>.\label{eq:forcing_coefficient}\\
a_{i} &= a_i(t).
\end{align}
\end{subequations}
\noindent Here the angle brackets denote the inner-product $\big<f,g\big> = \int_0^Lfg \text{ d}z$. Equation \eqref{eq:secondorderlinear} is a second order linear ODE with modal amplitudes $\bold{a}$, mass $\bold{M}$, damping $\bold{C}$ and stiffness $\bold{K}$ matrices and a time-dependent force vector $\ddot{v}_0\bold{d}$. When such a dynamical system is driven near or at its lowest resonance (the frequency corresponding to the smallest eigenvalue of the matrix $\bold{M}^{-1}\bold{K}$), quantitatively good results can be obtained even with a single-degree-of-freedom model, because all higher-order modes are expected to have negligible amplitudes in comparison. In addition, all higher-order modes are increasingly damped under the assumption that the damping is stiffness proportional, cf.\ Equations \eqref{eq:parameters} where $\bold{C}\propto \bold{K}$. In the contact phase the mechanical stop moves the imaginary part of the eigenspectrum to higher values, i.e.\ the stiffness increases. Combining \eqref{eq:stiffness_coefficient} with the natural assumption that the beam shape in the contact phase has larger curvature in the neighborhood of the mechanical stop, it can be expected that the lowest eigenfrequency of the beam increases when in contact.\par 
\noindent To obtain a single-degree-of-freedom model, we let $\eta(z,t)$ be the continuous, piecewise-defined, function
\begin{align}
\eta(z,t) =  \left\{
\begin{array}{c l l}
&a(t)\phi_f(z), & |a| < \Delta, \\
&\Delta\sgn(a)\phi_f(z)+(a(t)-\Delta\sgn(a))\phi_c(z) , & |a| \geq \Delta,
\end{array}
\right.
\end{align}
where $\sgn(\cdot)$ denotes the sign-function, subscripts $f$ and $c$ denote "free" and "contact", respectively, and $\Delta$ is the displacement amplitude of the mass point at which the beam grazes the mechanical stop. In what follows it is assumed that the cantilever displacement change is small compared to the change in curvature along the beam between solutions in free flight and contact phase. As a consequence we have
\begin{align}
0<\frac{|\eta_c-\eta_f|}{L^2(\eta_c''-\eta_f'')}\ll1,\label{eq:asymp_disp_curv}
\end{align}
where $\eta_f$ is a free configuration and $\eta_c$ is a contact configuration. Here, we assume that the beam mass is negligible relative to the lumped mass, i.e.,\ $\frac{\rho A_0 L}{m}\ll1$. Moreover, the following relations between the parameters in free flight and the contact phase are assumed:
\begin{subequations}
\begin{align}
c_c&=c_f(1+\beta), \\
k_c&=k_f(1+\alpha),\label{k2}\\
m_c&=m_f=m(1+\gamma), \\
d_c&=d_f(1+\nu),
\end{align}
\end{subequations}
where the relative changes in damping $\beta$ and stiffness $\alpha$ are large compared to $\gamma$ and $\nu$, respectively. Here we neglect the mass increment with $\gamma$ while retaining the forcing increment with $\nu$ for bookkeeping and later use. Dropping the subscripts, the system is then given by
\begin{subequations}
\label{dim_SDoF}
\begin{align}
&m\ddot{a}+c\dot{a}+ka+d\ddot{v}_0=0, &\text{ for }|a|<\Delta, \label{dim_SDoF_f}\\ 
&m\ddot{a}+c(1+\beta)\dot{a}+k(1+\alpha)a+\Gamma\sgn(a)+d(1+\nu)\ddot{v}_0=0,&\text{ for }|a|\geq\Delta. \label{dim_SDoF_c}
\end{align}
\end{subequations}
Here, the constant $\Gamma=-\alpha\Delta$ balances the static reaction forces of the beam at the mechanical stop, i.e.,\ the beam is bent and the elastic forces push it back towards the unbent state; this implies that the restoring force is continuous and piecewise linear. At $a=\pm\Delta$ the mechanical stops are reached. Let the clamped basepoint be harmonically moved, i.e.,\ $v_0=A\cos(\omega t)$. Rescaling and nondimensionalising the system by $t\mapsto \omega_1^{-1} t$ and $a\mapsto \Delta a$, where $\omega_1=\sqrt{\frac{k}{m}}$, and defining $\Omega=\frac{\omega}{\omega_1}$, $2\xi=\frac{c}{m\omega_1}$, and $I=\frac{d}{m}\frac{A}{\Delta}$, Equations \eqref{dim_SDoF} transforms into
\begin{subequations}
\label{eq:nondim_SDoF}
\begin{align}
&\ddot{a}+2\xi\dot{a}+a-I\Omega^2\cos(\Omega t)=0,&\text{ for }|a|<1, \label{eq:nondim_SDoF_f}\\
&\ddot{a}+2\xi(1+\beta)\dot{a}+(1+\alpha)a-\hspace{-0.1cm}\alpha\sgn(a)-\hspace{-0.1cm}I(1+\nu)\Omega^2\cos(\Omega t)=0,&\text{ for }|a|\geq1,\label{eq:nondim_SDoF_c}
\end{align}
\end{subequations}
which is the dynamical system analysed in Section \ref{sec:num_bif}.\par
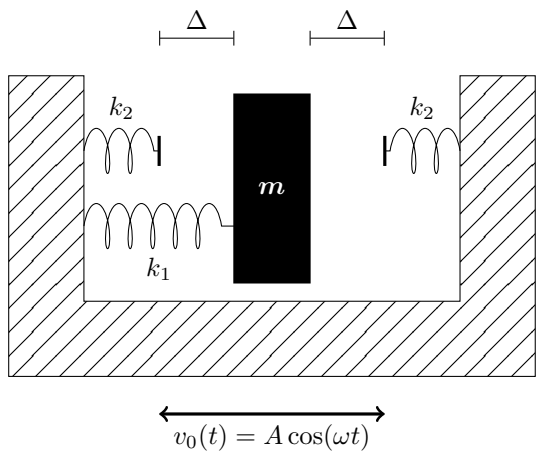
\begin{figure}[h]
\centering
\begin{tikzpicture}
\pgfdeclarepatternformonly{north east lines wide}%
   {\pgfqpoint{-1pt}{-1pt}}%
   {\pgfqpoint{10pt}{10pt}}%
   {\pgfqpoint{9pt}{9pt}}%
   {
     \pgfsetlinewidth{0.4pt}
     \pgfpathmoveto{\pgfqpoint{0pt}{0pt}}
     \pgfpathlineto{\pgfqpoint{9.1pt}{9.1pt}}
     \pgfusepath{stroke}
    }
\node (rect) at (3.5,2.5) [draw,thick,minimum width=1cm,minimum height=2.5cm,fill=black] {\textcolor{white}{{$\boldsymbol{m}$}}}; 
\draw[decoration={aspect=0.3, segment length=3mm, amplitude=3mm,coil},decorate] (1,2) --node[below,yshift=-3mm]{$k_1$} (3,2); 
\draw[decoration={aspect=0.3, segment length=3mm, amplitude=3mm,coil},decorate] (1,3) --node[above,yshift=3mm]{$k_2$} (2,3); 
\draw [very thick](2,2.8) -- (2,3.2);
\draw[decoration={aspect=0.3, segment length=3mm, amplitude=3mm,coil,mirror},decorate] (6,3) --node[above,yshift=3mm]{$k_2$} (5,3); 
\draw [very thick](5,2.8) -- (5,3.2);
\draw [pattern = north east lines wide] (0,0) -- (7,0) -- (7,4) -- (6,4) -- (6,1) -- (1,1) -- (1,4) -- (0,4) -- cycle;
\draw[<->,very thick] (2,-0.5)--node[below]{$v_0(t)=A\cos(\omega t)$}(5,-0.5);
\draw [|-|](2,4.5) -- node[above]{$\Delta$} (3,4.5);
\draw [|-|](4,4.5) -- node[above]{$\Delta$} (5,4.5);
\end{tikzpicture}
\caption{The mass/spring-system that is mechanically equivalent to the cantilever beam system in Figure \ref{fig:system_sketch}.}
\label{fig:free_body}
\end{figure}

An alternative approach yields an almost equivalent dynamical system: the model is derived by recognizing, a priori, that the mechanical system is analogous to that given in Figure \ref{fig:free_body}. Using Newton's second law for the mass, one arrives at the dynamical system
\begin{subequations}
\label{dim_SDoF_body}
\begin{align}
&m\ddot{a}+k_1a+\ddot{v}_0=0, &\text{ for }|a|<\Delta, \label{dim_SDoF__body_f}\\ 
&m\ddot{a}+k_2a+(k_1-k_2)\Delta\sgn(a)+m\ddot{v}_0=0,&\text{ for }|a|\geq\Delta. \label{dim_SDoF_body_c}
\end{align}
\end{subequations}
By adding damping, this mechanical analogue can be a good approximation to \eqref{dim_SDoF}, but in certain cases the changing geometry might be important; e.g.\ the modal mass, and forcing amplitude can vary. Furthermore, when using this direct ODE modeling approach it is difficult to assign a proper physical meaning to the mass $m$, and model and material parameters as well as geometrical parameters are not directly available. From a mathematical perspective, the derivation presented here and those given in, e.g.,\cite{MOON1983} and \cite{SHAW1985}, clarify underlying arguments that justify the reduction from PDE to low-dimensional ODEs.\par
\begin{figure}[h]
\centering
\subfloat[Subfigure 1 list of figures text][]{
\includegraphics[height=0.2\textwidth]{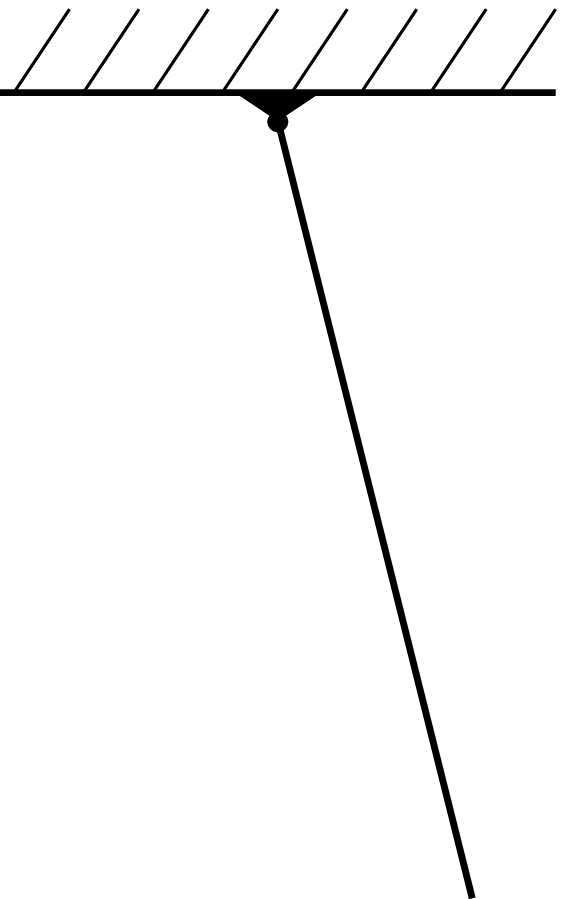}
\label{fig:mode01}}
\qquad
\subfloat[Subfigure 2 list of figures text][]{
\includegraphics[height=0.2\textwidth]{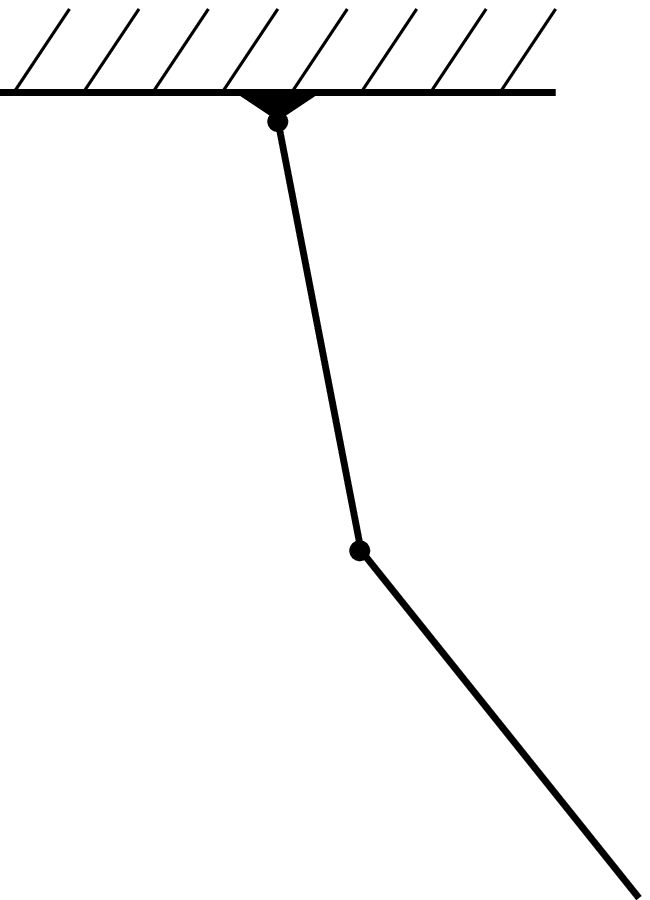}
\label{fig:mode02}}
\caption{The assumed configuration of the cantilever: Panel \protect\subref{fig:mode01} shows the one-element configuration in free flight; while Panel \protect\subref{fig:mode02} shows the two-element configuration in the contact phase.}
\label{fig:shapes}
\end{figure}

The model \eqref{eq:nondim_SDoF} was derived without any explicit use of specific shape functions nor numerical methods. To keep the derivation simple, we estimate the lowest-resonance frequencies in free flight and in contact in the next section. When comparing the theoretical findings to experimental data, a beam configuration as shown in Figure \ref{fig:shapes} is used.

\subsection{Parameter estimation}\label{subsec:Parameter_estimation}

\begin{figure}[t!]
\centering
\begin{picture}(7,5)
\put(0,0.25){ \includegraphics[width=0.5\textwidth]{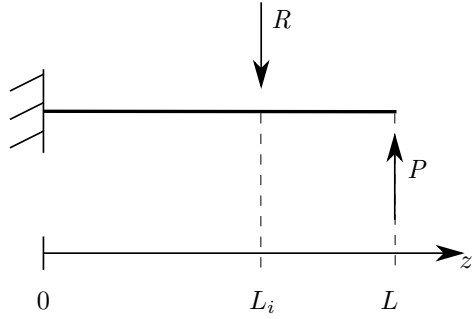}}
\put(0.48,-0.25){$0$}
\put(3.3,-0.25){$L_i$}
\put(5.05,-0.25){$L$}
\put(6.1,0.3){$z$}
\put(5.4,1.5){\textit{P}}
\put(3.6,3.5){\textit{R}}
\end{picture}
\caption{Cantilever beam exposed to transverse point forces $R$ and $P$. The reaction force $R$ acts as the mechanical stop, and $P$ is the point force acting on the mass point.}
\label{fig:static_beam_sketch}
\end{figure}

In free flight, the beam stiffness can be approximated by standard methods: The beam is driven at a frequency close to the first eigenfrequency; in this case the stiffness can be found via the lowest eigenfrequency of system \eqref{eq:PDE_comoving_all} without damping and forcing. We observe that for a single-mode expansion the undamped squared eigenfrequency is approximated via
\begin{align}
k/m=\frac{EI \int_0^L(\phi'')^2 \text{d}z}{M \phi^2(L_M)+\rho A\int_0^L\phi^2 \text{d}z}\, ,
\end{align}
where $k/m$ is calculated as $K_{11}/M_{11}$ from \eqref{eq:parameters}. Let the beam shape $\phi(z)$ be approximated by the cubic polynomial that solves the static problem, as is illustrated in Figure \ref{fig:static_beam_sketch} without the reaction force \textit{R}. The resulting the eigenfrequency is given by 
\begin{align}
f=\frac{1}{2\pi}\sqrt{\frac{k}{m}}\approx 8.4 \text{ Hz}.
\end{align}
This frequency deviates by about $10\%$ from the experimentally measured frequency of $7.75\text{ Hz}$. The estimation method is equivalent to the Rayleigh-Ritz method see, e.g., \cite{Courant1953} or, for a more exact approach, \cite{Meirovitch2001}; the gain in precision is $<1 \%$.\par
For vibrations involving contact, the estimation is more difficult due to the amplitude-dependent constraint which results in a nonlinear PDE. However, by assuming that this change in stiffness can be estimated by the static problem through the change in tip stiffness, the problem is linear. The equation that has to be solved is the differential equation for a slender bending beam under the assumption of small displacements (Euler-Bernoulli) and point loads, that is,
\begin{subequations}
\label{static_problem}
\begin{align}
&EIu''''=0 ,&\text{Euler-Bernoulli beam,}\\
&u(0)=u'(0)=0 ,&\text{clamped end,}\\
&u''(L)=0 ,&\text{moment free end,}\\
&EIu'''(L_{i})=R, \quad EIu'''(L)=P ,&\text{transverse point forces}.
\end{align}
\end{subequations}
This linear structural problem is statically indeterminate. The solution can be written as see, e.g., \cite{Timoschenko1997}
\begin{subequations}
\begin{align}
& u(\zeta)=u_1(\zeta)+u_2(\zeta), \text{ where }\zeta=\frac{z}{L}, \gamma = \frac{L_i}{L}\text{ and } \\
&u_1(\zeta) = \left\{
\begin{array}{c l l}
& -\frac{RL^3}{6EI}\zeta^2(3\gamma-\zeta), & \zeta \leq \gamma, \\
& -\frac{RL^3}{6EI}\gamma^2(3\zeta-\gamma), & 1 \geq \zeta > \gamma,
\end{array}
\right.\\
&u_2(\zeta) = \frac{PL^3}{6EI}\zeta^2(3-\zeta), \hspace{1.6cm} \zeta \leq 1.
\end{align}
\end{subequations}
\addtolength{\abovecaptionskip}{3mm}
The relative change in tip stiffness is given by $\alpha$, where tip stiffness is defined as 
\begin{align*}
k=\left.\pdd{P}{u}\right|_{\zeta=1}.
\end{align*}
When the beam does not touch the mechanical stop the reaction force is zero, i.e.,\ $R=0$ and, therefore, the tip stiffness is given by
\begin{align}
k_1=\left.\pdd{P}{u}\right|_{\zeta=1}=\frac{3EI}{L^3} \label{k1_static}
\end{align}
for free flight. When the beam touches the mechanical stop one has to know $R=R(P)$. At the impact,  $u(\gamma)=\delta$ and this implies a linear relationship between $P$ and $R$ given by
\begin{align*}
u(\gamma)=\delta= -\frac{RL^3}{3EI}\gamma^3+ \frac{PL^3}{6EI}\gamma^2(3-\gamma).
\end{align*}
From this expression it follows that the tip stiffness in the contact phase is
\begin{align}
k_2=\left.\pdd{P}{u}\right|_{\zeta=1}=\frac{3EI}{L^3} \left( 1-\frac{1}{4}\left(3-\gamma\right)^2\gamma\right)^{-1}. \label{k2_static}
\end{align}
Combining \eqref{k2}, \eqref{k1_static} and \eqref{k2_static} leads to
\begin{align}
\alpha = \left( 1-\frac{1}{4}\left(3-\gamma\right)^2\gamma\right)^{-1}-1 =\left( 1-\frac{1}{4}\left(3-\frac{L_i}{L}\right)^2\frac{L_i}{L}\right)^{-1}-1.\label{stiffness}
\end{align}
It appears that $\alpha \rightarrow 0$ for $L_i \rightarrow 0$, and $\alpha \rightarrow \infty$ for $L_i \rightarrow L $.
In other words, the change in stiffness vanishes as the mechanical stop is moved towards the clamped end, and it is unbounded as the mechanical stop approaches the free tip. This agrees with the physics of the problem. When the mechanical stop is at the clamped end, it is not going to impact with the beam so the change in tip stiffness vanishes. When the mechanical stop approaches the tip of the beam, it will become harder to move the tip of the beam transversely, to the point where it is not possible to do so as the mechanical stop is blocking any motion of the tip.\par
Since the damping is chosen  proportional to the stiffness \eqref{eq:PDE_comoving_all}, the change in damping has the same order of magnitude as $\alpha$. In Table \ref{tab:Parameter_estimation} we see that $\alpha$ deviates by about $20\%$ from the experimentally fitted value.\par
 \begin{table}[htdp]
\begin{center}
\begin{tabular}{|l|c|c|}
\hline & Theory & Experiment \\ \hline
$f$& $\sim$ 8.4 Hz & 7.75 Hz \\ \hline
$\alpha$ & 4.9 & 5.9 \\ \hline
$\xi$ & --- &  3 \% \\ \hline
$\beta$ & $\mathcal{O}(1)$ &  0.885 \\ \hline
\end{tabular}
\end{center}
\vspace{-0.7cm}
\caption{Stiffness parameters $f$ and $\alpha$, and damping parameters $\xi$ and $\beta$. The parameter $f$ is the natural frequency of the cantilever. Experimental values are model-fitted values.}
\label{tab:Parameter_estimation}
\end{table}
In summary, the mathematical model \eqref{eq:nondim_SDoF} has been derived for a clamped beam model with symmetric mechanical stops. The parameters of the model can be split into two main groups as follows:
\begin{itemize}
\item Fixed parameters: damping $\xi$ and ratio $\beta$ of damping constants are obtained from experimental data.
\item Free parameters: forcing amplitude $I$, frequency $\Omega$, stiffness ratio $\alpha$, and relative forcing amplitude $\nu$.
\end{itemize}
The damping $\xi$ and the natural frequency $f$ (Table \ref{tab:Parameter_estimation}) in free flight are fitted to the experimental data. The relative stiffness $\alpha$ is considered free, but in the actual experiment it was not changed. 

\section{Smoothing of the piecewise-smooth impact model}\label{sec:smoothing}

\noindent In the present case of a cantilever beam with mechanical stops a smooth representation is constructed and analyzed. The smoothing is realized through a nonlinear homotopy or, equivalently, a smooth switching function. To this end, we consider a differentiable nonlinear scalar switching function $\text{H}=\text{H}(\bold{x},p)$ where H converges pointwise to a Heaviside-type step function as $p\rightarrow \infty$. Using $\bold{x}=(x_1,x_2)^T$, where $x_1 = a$ and $x_2 = \dot{a}$, Equation \eqref{eq:nondim_SDoF} is reformulated as the first-order system
\begin{align}
\dot{\bold{x}} = \text{f}(\bold{x},t,\bm{\lambda}) = \left\{ 
\begin{array}{c l}
& \text{f}_1(\bold{x},t,\bm{\lambda})),\quad  \quad | x_1 | \leq 1, \\
& \text{f}_2(\bold{x},t,\bm{\lambda})),\quad  \quad | x_1 | > 1,
\end{array}
\right.
\end{align}
where $\text{f}:\mathbb{R}^{2}\times\mathbb{R}\times \mathbb{R}^{k+1}\mapsto \mathbb{R}^{2}$ is piecewise smooth. 
\noindent The smooth dynamical equivalent is obtained by the following homotopy
\begin{align}
\label{nonlinear_homotopy}
\dot{\bold{x}} = \text{g}(\bold{x},t,\bm{\lambda},p) := \text{H}\text{f}_1+(1-\text{H})\text{f}_2,
\end{align}

\noindent so that $\text{g}:\mathbb{R}^{2}\times\mathbb{R}\times \mathbb{R}^{k+1}\times \mathbb{R}\mapsto \mathbb{R}^{2}$ is now a smooth function. This approach has been taken in \cite{Szalai2009a}; see also \cite{Kuznetsov2003}.The choice of $\text{H}(\bold{x},p)$ is a modeling decision and, in general, one can smooth a given sequence of impacts in a piecewise-defined system by smooth functions that approximate the nonsmooth right-hand sides. Well-known examples are the Bump function (compact support) and the hyperbolic tangent function. In general, a piecewise-smoothly defined system has smooth representations given by sums, $\text{g}(\bold{x},\bm{\lambda},\bold{p})=\sum\limits_i \text{H}_i(\bold{x},\bm{\lambda},p_i) \text{f}_i(\bold{x},\bm{\lambda})$. Note that the smooth switching functions can also depend on the system parameters and in some cases this might be of interest if, e.g.,\ the natural switching occurs over a predefined range depending on a model-specific geometric parameter. If the smoothing procedure is performed in order to construct a smooth equivalent to the nonsmooth system then a necessary condition is that there exists smoothing homotopy parameter values $p_i^*$ such that the bifurcation diagrams are robust for all $p_i>p_i^*$.\par
Smoothing processes introduce additional parameters $p_i$ that may effect the solution and bifurcation structure. Hence, it is important to analyse, e.g., how special points depend on the smoothing homotopy parameters and whether they are robust with respect to changes in the parameters $p_i$. In the present case, with one smoothing homotopy parameter $p$, we are interested in quantitative measures, such as the displacement of the cantilever, which should not change for sufficiently large $p$. From the point of view of numerical computations it is not tractable to make $p$ extremely large, because this introduces very large derivatives near the nonsmooth events. Therefore, the strategy is to find a suitable value of $p$ that does not induce artificial and unwanted dynamics or skew the quantitative results in the parameter domains of interest.\par
\begin{figure}[h!]
\centering
\includegraphics[width=0.7\textwidth]{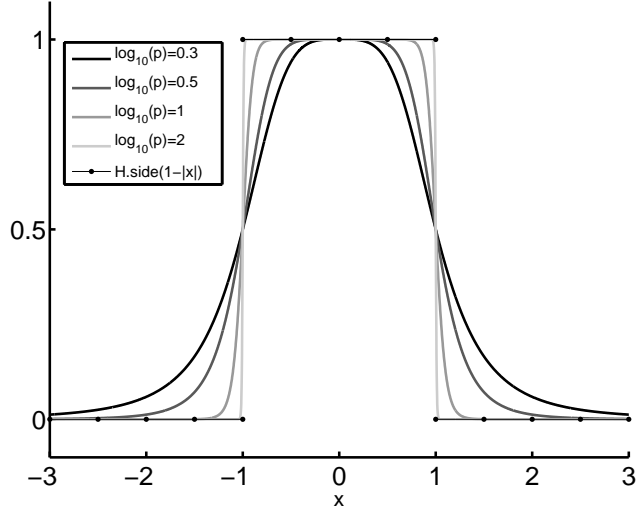}
\caption{Smooth switching functions $\text{H}(x,p)$ for $\log_{10}(p)\in \{0.3,0.5,1,2\}$ compared to the Heaviside step function $\text{Heaviside}(1-|x|)$.}
\label{fig:homfunc}
\end{figure}
\begin{figure}[h]
\centering
\includegraphics[width=0.7\textwidth]{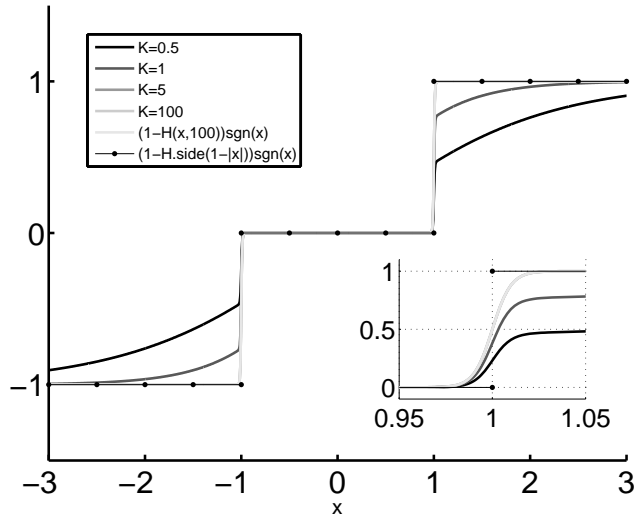}
\caption{Smooth switching term $(1-\text{H}(x,p))\tanh(Kx)$ for $p=100$ and $K\in \{0.5,1,5,100\}$ compared with piecewise-smooth switching functions $(1-\text{H}(x,100))\sgn(x)$ and $(1-\text{Heaviside}(1-|x|))\sgn(x)$.}
\label{fig:homtanh}
\end{figure}

We consider the dynamical behavior of the smoothed equivalent of system \eqref{eq:nondim_SDoF} with parameters $\bm{\lambda}=(\xi,\beta,\alpha,\Omega,I,\nu)$. Specifically, we consider the cantilever beam with stops
\begin{align}
\label{eq:BD_system}
\dot{\vt{x}} = \text{H}\text{f}_1+(1-\text{H})\text{f}_2,
\end{align}
with
\begin{align}
\text{f}_1(\bold{x},t,\bm{\lambda})&=\begin{pmatrix} 0 &1\\ -1 &-2\xi \end{pmatrix} \bold{x}+I\Omega^2 \cos(\Omega t) \begin{pmatrix}  0 \\ 1 \end{pmatrix}, \nonumber \\
\text{f}_2(\bold{x},t,\bm{\lambda})&=\begin{pmatrix} 0 &1\\ -(1+\alpha) &-2\xi(1+\beta) \end{pmatrix} \bold{x}+\alpha \begin{pmatrix} 0 \\ \tanh(Kx_1) \end{pmatrix}+I\Omega^2 \cos(\Omega t) \begin{pmatrix}  0 \\ 1+\nu \end{pmatrix}, \nonumber
\end{align}
and
\begin{align}
\text{H}(\bold{x},p) &= \Big(1+(x_1^2)^{p}\Big)^{-1}. \nonumber
\end{align}
Our choice of the smooth switching function $\text{H}$ is visualized in Figure~\ref{fig:homfunc} for several values of $p$. In the present case, the smoothing depends only on the displacement of the beam and the smoothing homotopy parameter $p$.
Notice that the definition of $\text{f}_2(\bold{x},t,\bm{\lambda})$ is modified slightly by approximating the discontinuous sign-function in Equation \eqref{eq:nondim_SDoF_c} with the smooth hyperbolic tangent function. This is a necessary step in order to make the dynamical system smooth. In the nonimpacting region any term in $\text{f}_2(\bold{x},t,\bm{\lambda})$ is multiplied by the smooth switching term $1-\text{H}(\bold{x},p)$ which converges to zero pointwise as $p\rightarrow\infty$. Therefore, it only matters how well the hyperbolic tangent function approximates the sign-function for $|x|\geq1$. In Figure \ref{fig:homtanh}, the piecewise-constant switching term and $(1-\text{H}(x,p))\sgn(Kx)$ is plotted together with the smooth approximations $(1-\text{H}(x,p))\tanh(Kx)$ with $\log_{10}(p)=2$ and $K=0.5$, $K=1$, $K=5$ and $K=100$.
The smooth approximations are visually indistinguishable for $K=5$ and $K=100$; we choose to fix $K=100$. In the impacting regions the error between the sign-function and the hyperbolic tangent function with $K=100$ is \centerline{$1-\tanh(100|x_1|)<10^{-15}$ for $|x_1|>1$.} Since the local variations are below the numerical tolerances, we can conclude that further investigations of the dependence of the smoothing on the smoothing parameter $K$ is not necessary.

\section{Numerical bifurcation analysis}\label{sec:num_bif}

\begin{figure}[h]
\centering
\includegraphics[width=0.45\textwidth]{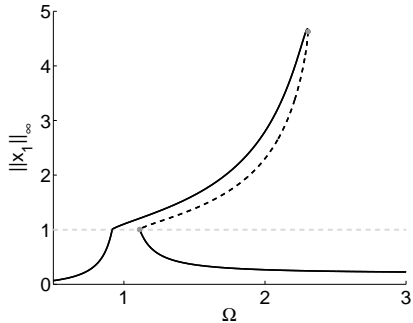}
\caption{Nonlinear frequency response diagram of system \eqref{eq:BD_system} with $\nu = 0$, $\alpha = 5.9$ and $\log_{10} (p) = 2$. Stable branches are solid, unstable branches are dashed and the two fold points are marked with grey bullets. The other parameters are $(\xi,\beta,I)=(0.03,0.885,0.2)$.}
\label{fig:alpha5nu0_eps_2_freqresp}
\end{figure}

In this section we investigate the dynamical behavior of the smooth system using the numerical continuation and bifurcation analysis package AUTO \cite{auto}. We continue periodic orbits in one and two parameters of the smooth system \eqref{eq:BD_system} at the 1:1 resonance. As mentioned earlier, the mechanical system is, a harmonically-driven single-degree-of-freedom nonlinear oscillator for which the nonlinearity is hardening, because the effective beam stiffness is larger in the contact phase. When such a system is driven with a large enough forcing amplitude near its natural frequency (eigenfrequency) the usual response is a nonlinear resonance peak. Figure \ref{fig:alpha5nu0_eps_2_freqresp} shows one of the nonlinear resonance peaks, where $\Omega$ denotes the scaled forcing frequency and $||x_1||_\infty$ denotes the maximal displacement amplitude of the periodic orbits. Stable and unstable branches are indicated by solid and dashed lines, respectively. For a range of frequencies, we have a bistable system with three coexisting periodic orbits, two stable orbits separated by an unstable orbit. There are two bifurcation points marked by the two grey bullets, which are fold points (limit point/saddle point/turning point). Together, the fold points organize a hysteresis loop that determines the possible behavior of the system when keeping all parameters fixed except for the forcing frequency. Starting with a relatively large forcing frequency, the system stabilizes to a solution on the lower branch. Continuously decreasing the forcing frequency causes a discontinuous (increase) jump in the response amplitudes when reaching the left-most fold point. Similarly, starting from a relatively small forcing frequency, and continuously increasing the driving frequency results in a discontinuous (decrease) jump in response amplitude when passing the right fold point on the upper stable branch.\par 

In Section \ref{subsec:BD_1} we investigate the behavior of the system without considering a change in the forcing amplitude at the mechanical stop, i.e.,\ for $\nu=0$. In Section \ref{subsec:BD_3} the case of doubling the forcing amplitude in the contact phase is considered, i.e.,\ we set $\nu=1$. In this case, the system has two additional fold points and the bifurcation diagram exhibits an isola of periodic solutions that is not connected to the main solution branch. In what follows we neglect the discontinuous damping terms when the regularity is stated, this is motivated by the fact that the damping is orders of magnitude smaller than all other terms. 

\subsection{Continuous forcing amplitude}\label{subsec:BD_1}

\begin{figure}[h]
\centering
\subfloat[Subfigure 1 list of figures text][]{
\includegraphics[width=0.45\textwidth]{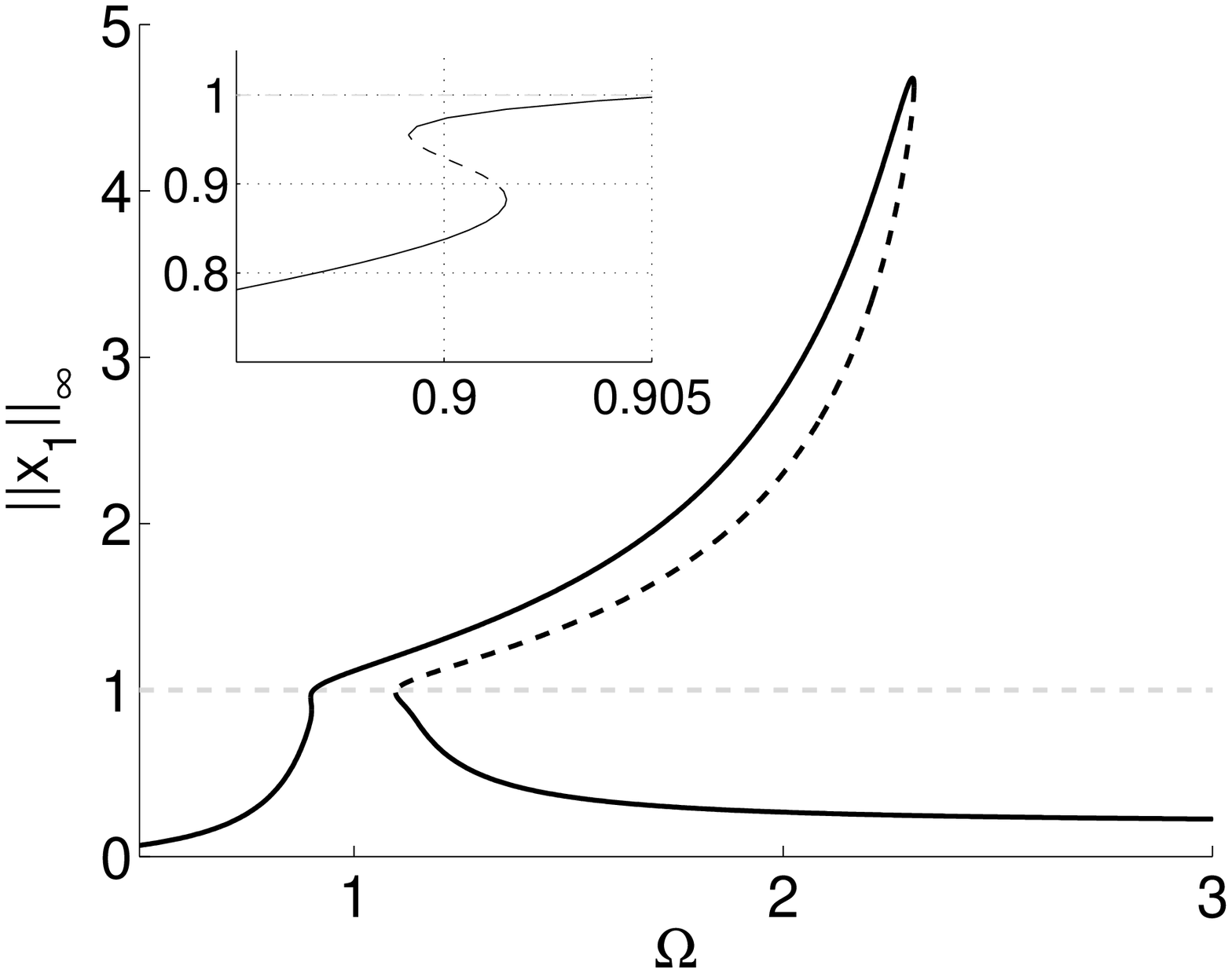}
\label{fig:alpha5nu0_freqresp}}
\hfill
\subfloat[Subfigure 2 list of figures text][]{
\includegraphics[width=0.45\textwidth]{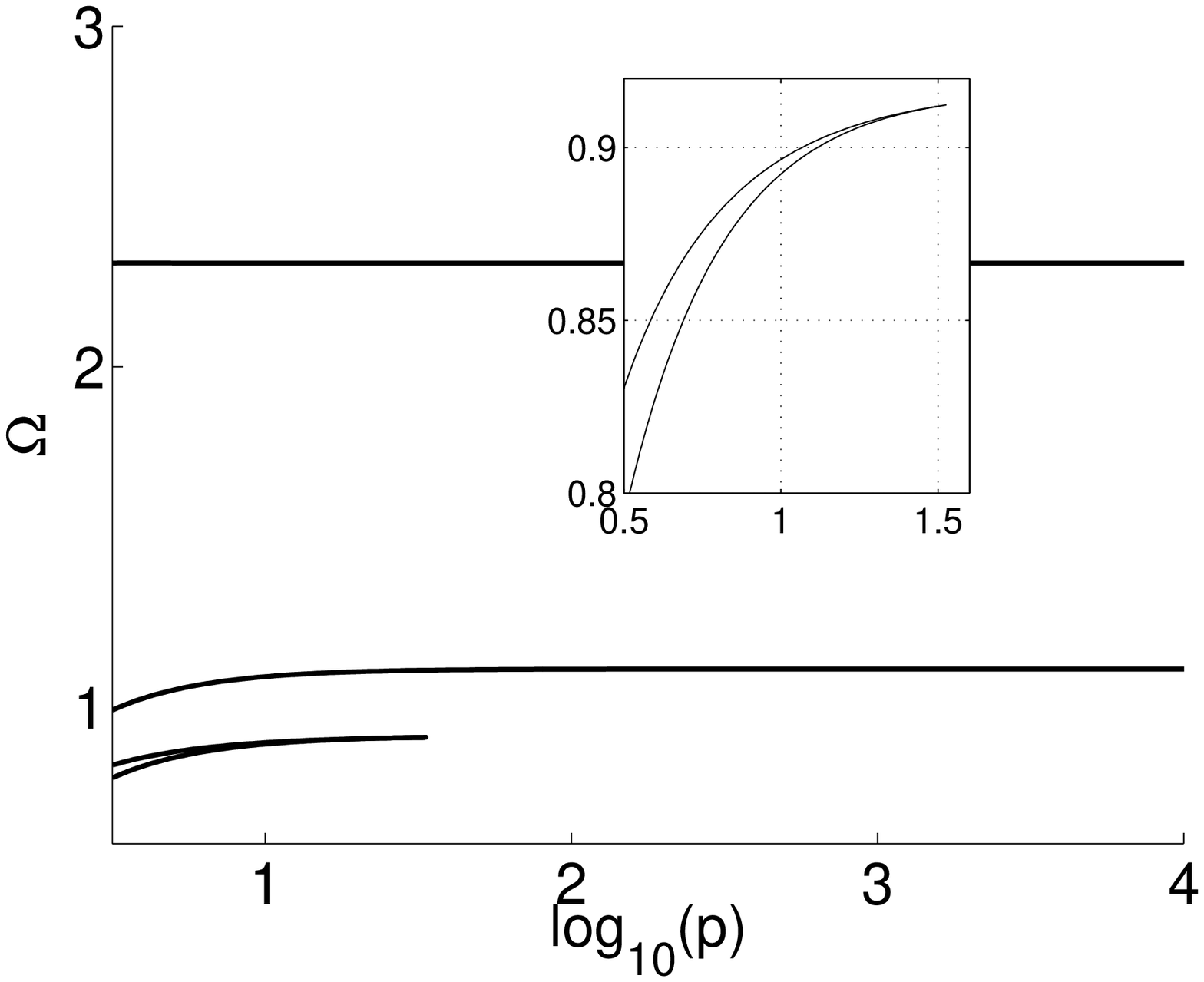}
\label{fig:alpha5nu0}}
\caption{System \eqref{eq:BD_system} with $\nu = 0$ and $\alpha = 5.9$; Panel \protect\subref{fig:alpha5nu0_freqresp} shows the frequency response diagram for $\log_{10} (p) \approx 1.1$; Panel \protect\subref{fig:alpha5nu0} shows the loci of fold points in the ($\log_{10}(p),\Omega$)-plane. The remaining parameter set is $(\xi,\beta,I)=(0.03,0.885,0.2)$.}
\label{fig:globfig01}
\end{figure}

First we consider the case of a continuous forcing amplitude, i.e.,\ $\nu = 0$, so that system  \eqref{eq:nondim_SDoF} is Lipschitz continuous. Moreover, we choose a moderately stiff impact using an intermediate value of $\alpha = 5.9$. Based on Equation \eqref{stiffness}, $\alpha=5.9$ geometrically corresponds to a placement of the mechanical stops at $3/5$-length of the beam (from the clamped end $z=0$). Figure \ref{fig:globfig01}\protect\subref{fig:alpha5nu0_freqresp} shows the corresponding one-parameter bifurcation diagram; as in Figure \ref{fig:alpha5nu0_eps_2_freqresp}, $||x_1||_\infty$ denotes the maximal amplitude of the periodic orbits and the bifurcation parameter is the scaled frequency $\Omega$. At first sight, the bifurcation diagram shows a typical nonlinear resonance like the one shown in Figure \ref{fig:alpha5nu0_eps_2_freqresp}. However, two additional very localized folds are robustly present on the upper branch; see the enlargement in Figure \ref{fig:globfig01}\protect\subref{fig:alpha5nu0_freqresp}. Since we are smoothing a piecewise-smooth system, we also investigate how the bifurcation diagram depends on the smoothing homotopy parameter $p$. To this end, we follow the loci of all four fold points in two parameters. Figure \ref{fig:globfig01}\protect\subref{fig:alpha5nu0} shows the bifurcation diagram projected onto the $(\log_{10}(p), \Omega)$-parameter plane. The two upper fold points are robust, in the sense that they are practically independent of $p$ as $p\rightarrow \infty$, while the two lower fold points disappear in a cusp for $\log_{10} (p) \approx 1.5$. Hence, for $p$ large enough, the smoothed system exhibits the same expected bifurcation structure as system \eqref{eq:nondim_SDoF}. We conclude that, with a Lipschitz-continuous vector field for $\nu=0$, the smoothed system shares the same bifurcation structure as the Duffing oscillator near resonance, because the two localized fold points on the left of the resonance are induced by the smoothing. In another setting, it may be advantageous to create this smoothing effect deliberately by design; however, for the mechanical system under consideration here, this is unwanted behavior.\par
\begin{figure}[h]
\centering
\subfloat[Subfigure 1 list of figures text][]{
\includegraphics[width=0.45\textwidth]{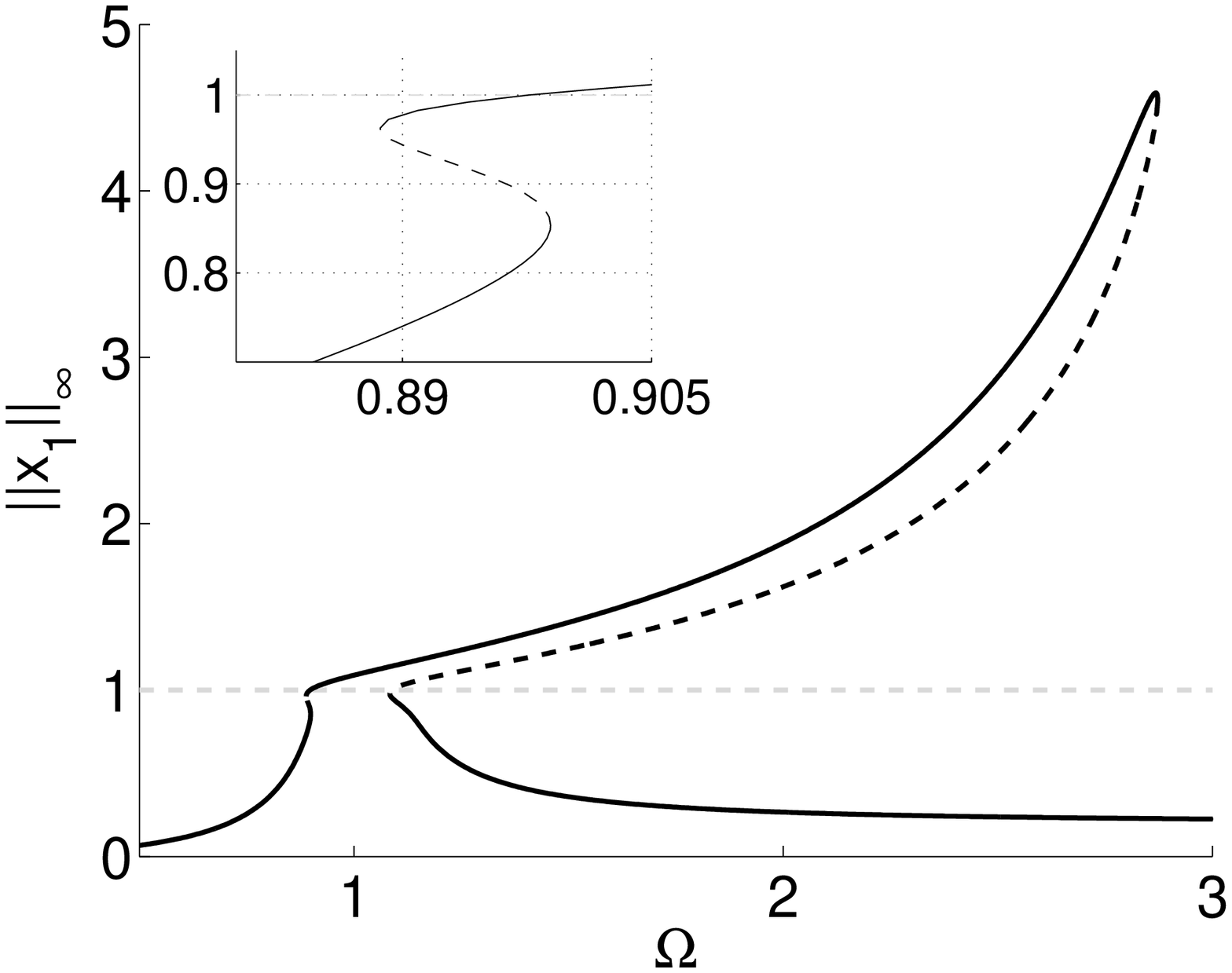}
\label{fig:alpha10nu0_freqresp}}
\hfill
\subfloat[Subfigure 2 list of figures text][]{
\includegraphics[width=0.45\textwidth]{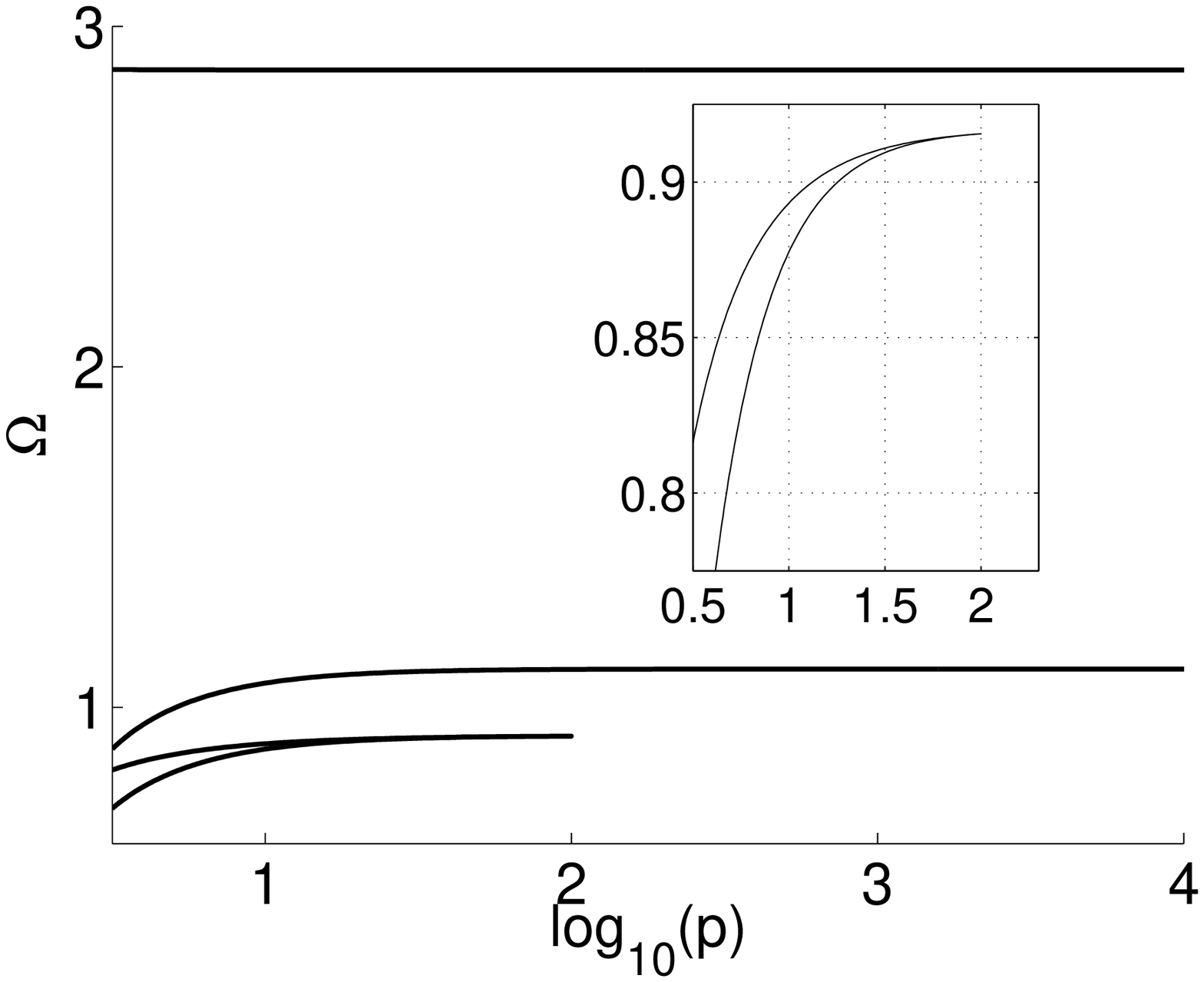}
\label{fig:alpha10nu0}}
\caption{System \eqref{eq:BD_system} with $\nu = 0$ and $\alpha = 10$; Panel \protect\subref{fig:alpha5nu0_freqresp} shows the frequency response diagram for $\log_{10}(p)\approx 1.1$; Panel \protect\subref{fig:alpha5nu0} shows the loci of fold points in the $\left(\log_{10}(p),\Omega\right)$-plane. The remaining parameter set is $(\xi,\beta,I)=(0.03,1.5,0.2)$.}
\label{fig:globfig02}
\end{figure}

For $\alpha = 10$, still with $\nu=0$, the change in stiffness is almost twice as large as for $\alpha = 5.9$. Based on Equation \eqref{stiffness}, $\alpha=10$ geometrically corresponds to a placement of the mechanical stops at $2/3$-length of the beam (from the clamped end $z=0$). There is no qualitative difference between Figures~\ref{fig:globfig01}~and~\ref{fig:globfig02}. Note from the enlargements that the two localized fold points are both below threshold of the mechanical stop; this indicates that they are induced by the smoothing, since the piecewise-smooth system is linear if the mechanical stops are not impacted upon, i.e., for solutions with $||x||_\infty<1$.\par
\begin{figure}[h]
\centering
\subfloat[Subfigure 1 list of figures text][]{
\includegraphics[width=0.45\textwidth]{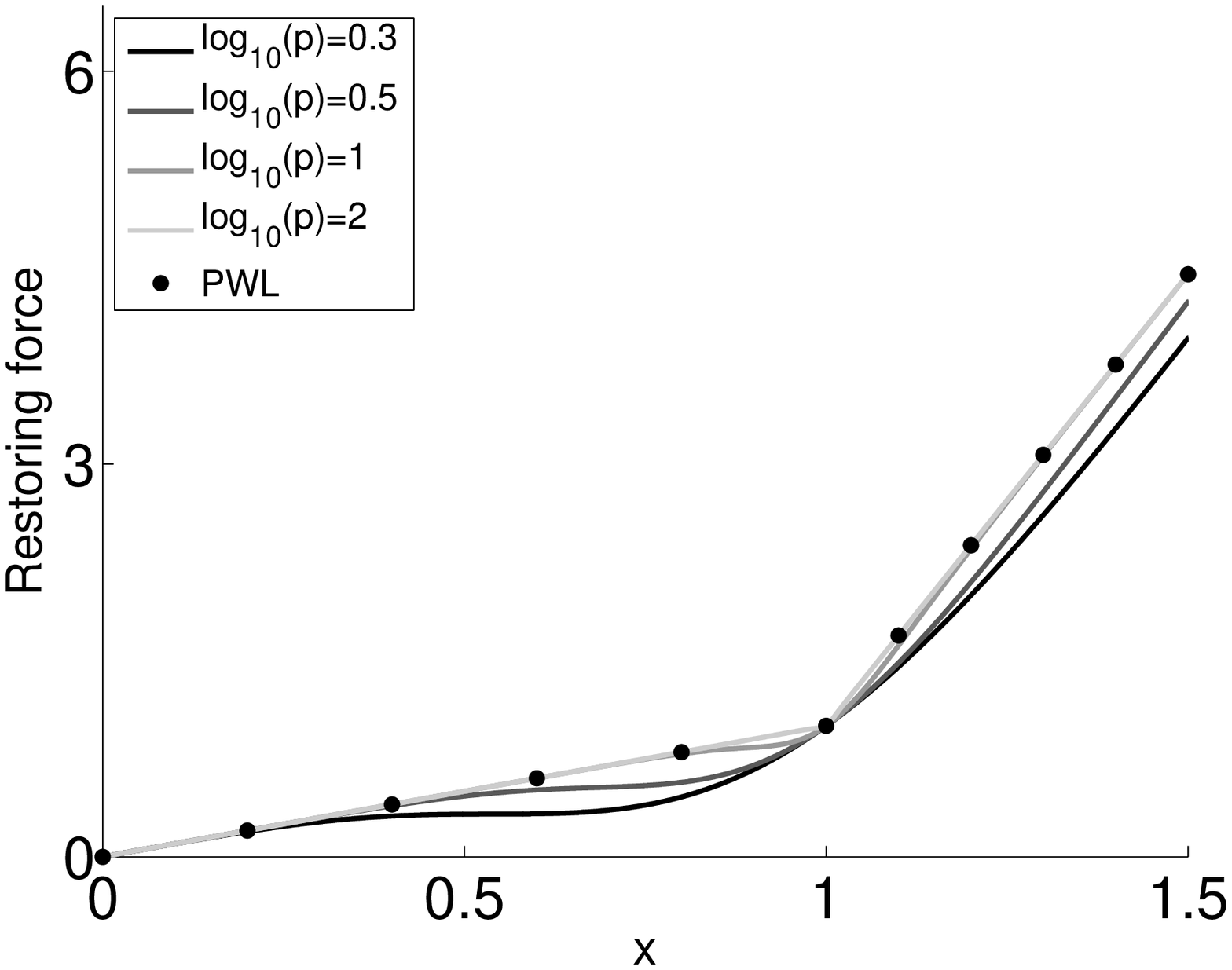}
\label{fig:force_vs_hom_alpha5_eps05}}
\hfill
\subfloat[Subfigure 2 list of figures text][]{
\includegraphics[width=0.45\textwidth]{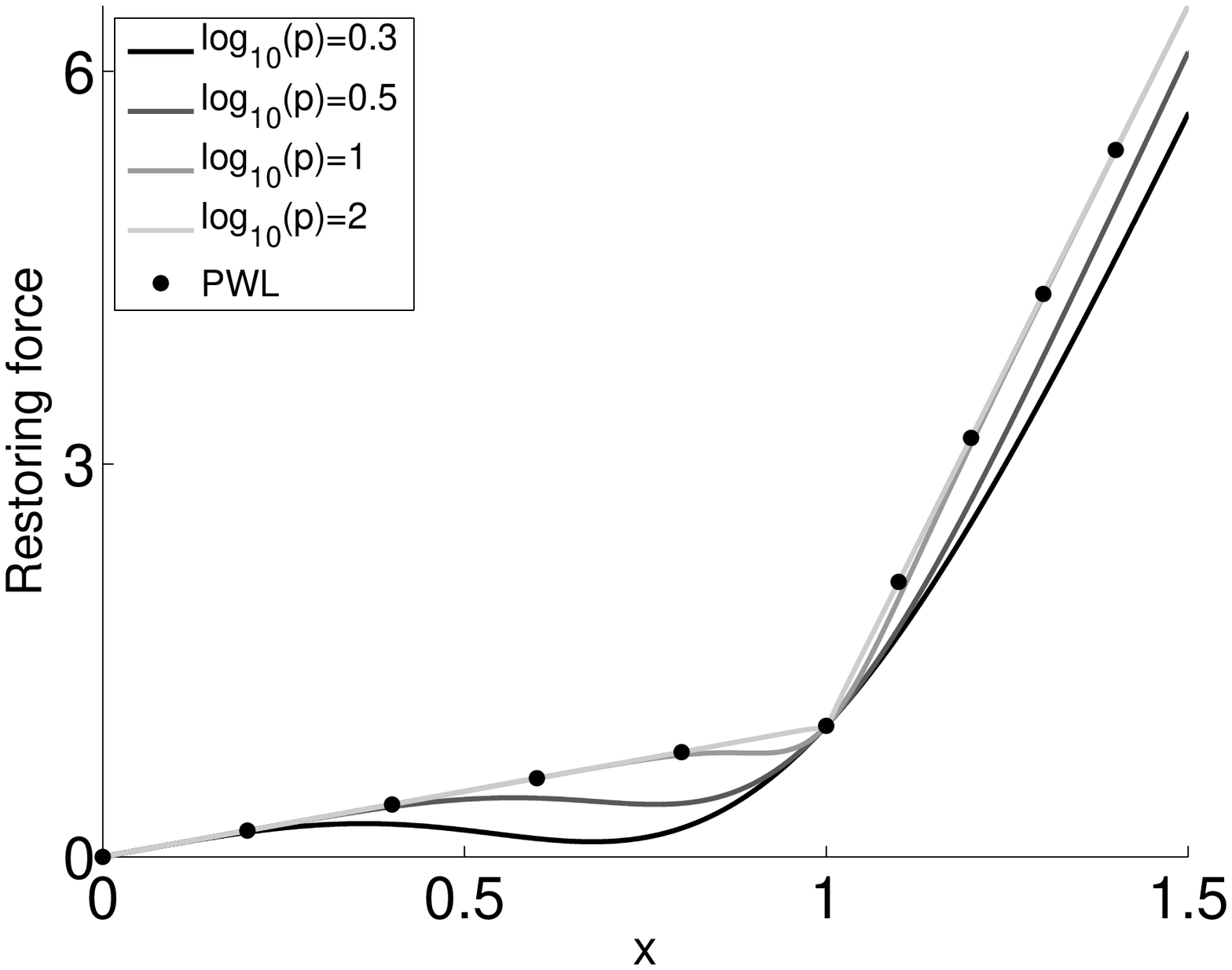}
\label{fig:force_vs_hom_alpha10_eps05}}
\caption{Elastic restoring force as a function of transverse displacement. Panel~\protect\subref{fig:force_vs_hom_alpha5_eps05} shows the moderate impact, $\alpha = 5.9$. Panel~\protect\subref{fig:force_vs_hom_alpha10_eps05} shows the hard impact $\alpha = 10$. The legends refer to different smoothing levels of the piecewise-linear force, i.e.,\ different $p$. PWL denotes the piecewise-linear force.}
\label{fig:globfig_force_vs_hom}
\end{figure}

The additional fold points emerge, because the smooth switching is not sufficiently localized. Consequently, the restoration force features a softening effect that is unphysical. For each given $\alpha$ there is a suitable choice of $p$, but as $\alpha$ is increased, the smoothing must be more localized, i.e.,\ $p$ must be increased. If $p$ is not sufficiently large the softening effect becomes more pronounced, and in the present case this is the sole mechanism for the observed smoothing-induced bifurcations; softening and hardening effects are explained, e.g.,\ in \cite{Hale1969} and \cite{Nayfeh1995}. Figure \ref{fig:globfig_force_vs_hom} illustrates how, for small $p$, the smooth approximations of the piecewise-linear elastic restoring force and/or its derivative dips (i.e.,\ spring softens) in a neighborhood of $x=1$. In Figure \ref{fig:alpha10nu0} the fold points are continued in $(\log_{10}(p),\Omega)$. It is observed that they disappear in a cusp for $\log_{10}(p)\approx 2$, hence, the smoothing should be chosen such that $\log_{10}(p)>2$.

\subsection{Discontinuous forcing amplitude}\label{subsec:BD_3}

\begin{figure}[h]
\centering
\subfloat[Subfigure 1 list of figures text][]{
\includegraphics[width=0.45\textwidth]{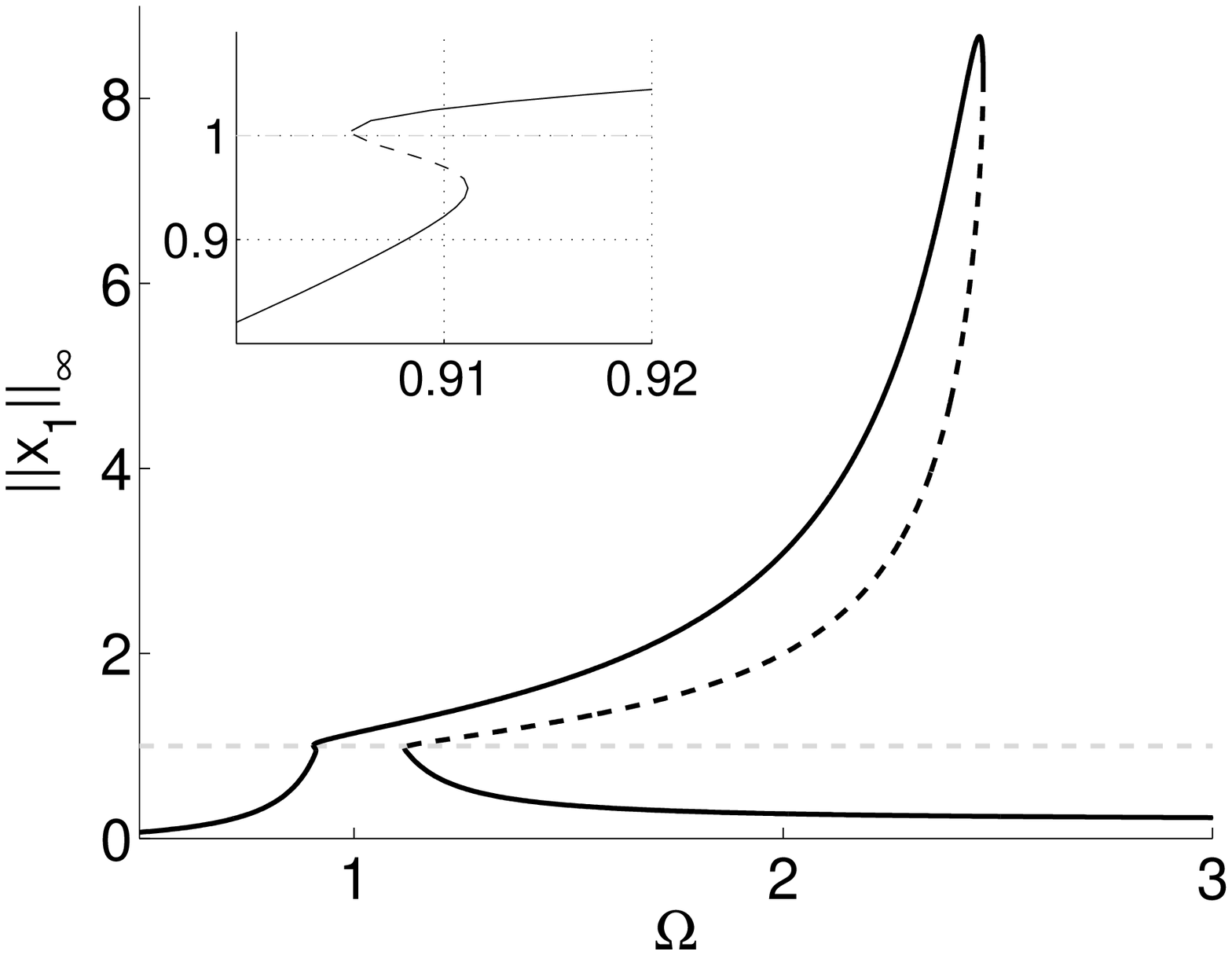}
\label{fig:alpha5nu1_freqresp}}
\hfill
\subfloat[Subfigure 2 list of figures text][]{
\includegraphics[width=0.45\textwidth]{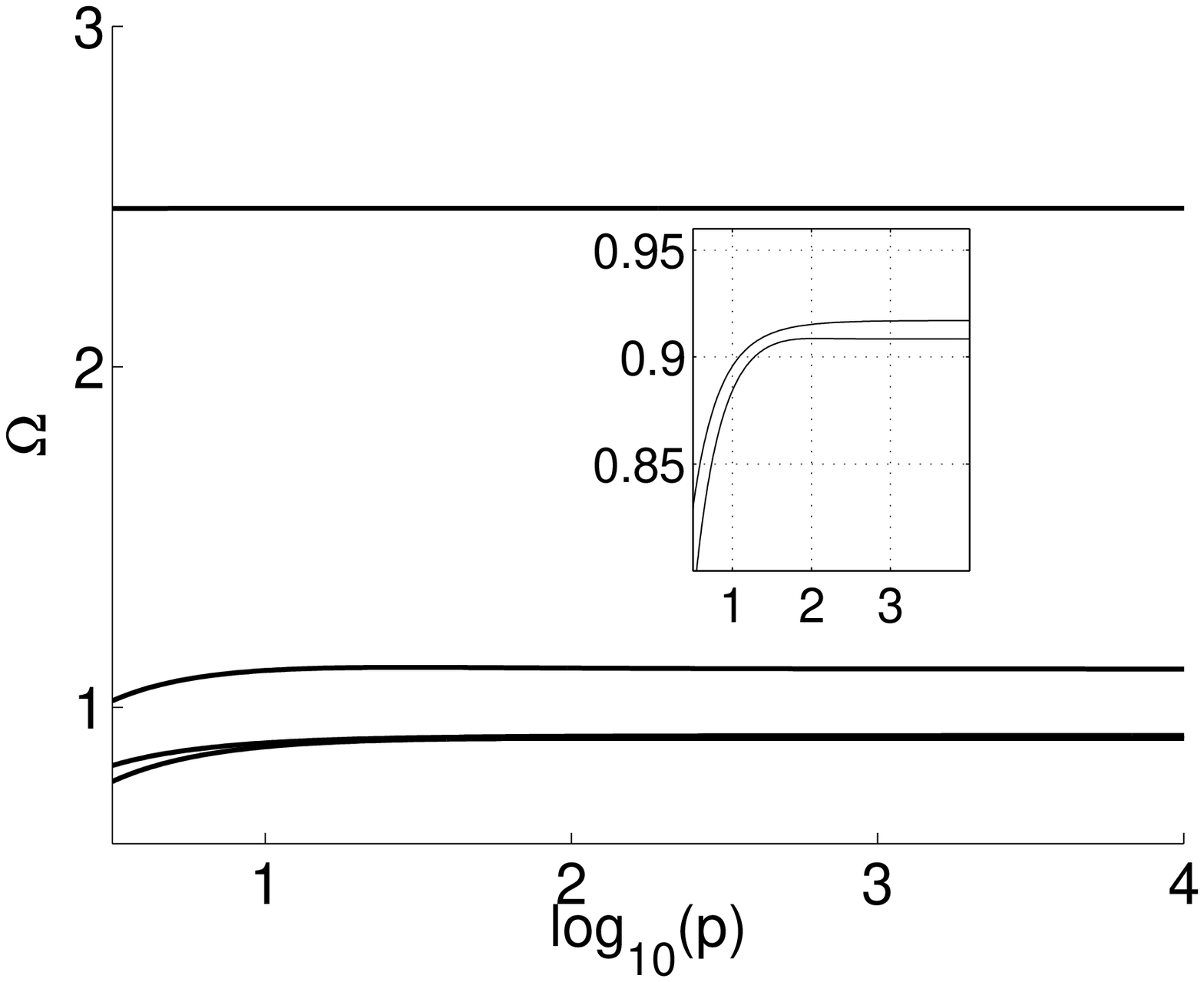}
\label{fig:alpha5nu1}}
\caption{System \eqref{eq:BD_system} with $\nu=1$ and $\alpha = 5.9$. Panel \protect\subref{fig:alpha5nu1_freqresp} shows the frequency response diagram for $\log_{10}(p)=1.5$. Panel \protect\subref{fig:alpha5nu1} shows the loci of fold points in the $(\log_{10}(p),\Omega)$-parameter plane; all four branches are disconnected. The remaining parameter set is $(\xi,\beta,I)=(0.03,0.885,0.2)$.}
\label{fig:globfig03}
\end{figure}

We next consider the case where the the forcing amplitude is such that system \eqref{eq:nondim_SDoF} is discontinuous, by setting $\nu = 1$. We first set $\alpha=5.9$, so that the cantilever is undergoing moderately stiff impacts. These parameter values will also be used when the system is compared with experimental data in Section \ref{sec:exp_comp}. Figure \ref{fig:globfig03}\protect\subref{fig:alpha5nu1_freqresp} shows that, for $\nu=1$, two fold points are again introduced in the same region as in Figure \ref{fig:globfig01}\protect\subref{fig:alpha5nu0_freqresp} and the qualitative structure is the same. However, these fold points are not induced by the smoothing process. This is illustrated by the two-parameter bifurcation diagram in Figure \ref{fig:globfig03}\protect\subref{fig:alpha5nu1}, which shows that the two additional fold points do not disappear in a cusp as $p$ is increased. For this particular parameter set there are four limit points for any $p$.\par
As was mentioned earlier, it may be of interest to include dependencies of model parameters in the smoothing functions. For example, the smooth switching functions generally depend on model parameters, i.e.,\ $p=\Theta(\lambda)$ where $\Theta$ is designed by a-priori model properties or from numerical continuation results. Our analysis shows that the smooth switching must be scaled according to the stiffness ratio $\alpha$, i.e.,\ $p=\Theta(\alpha)$; however, since the experimental data is for fixed $\alpha$, we need not construct $\Theta(\alpha)$ explicitly for more values of $\alpha$ than are already investigated. We mention here the modeling and analysis of the dynamics of the main landing gear on a plane \cite{C.Howcroft}, where such a smoothing dependence is analyzed further, because, in that case, the smooth switching function naturally depends on an active bifurcation parameter.\par

\begin{figure}[h]
\centering
\subfloat[Subfigure 1 list of figures text][]{
\includegraphics[width=0.45\textwidth]{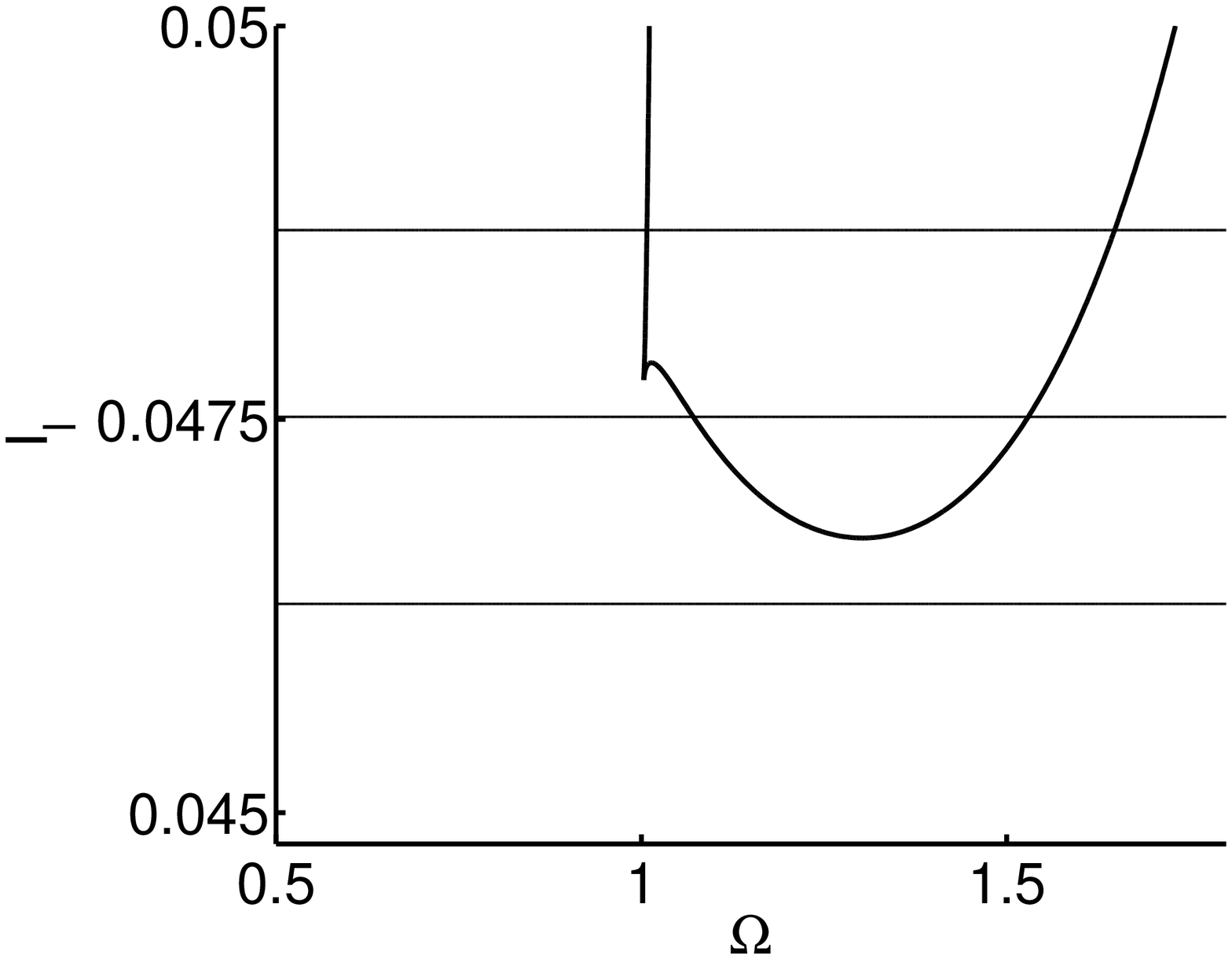}
\label{fig:lp_isola_curve}}
\hfill
\subfloat[Subfigure 2 list of figures text][]{
\includegraphics[width=0.45\textwidth]{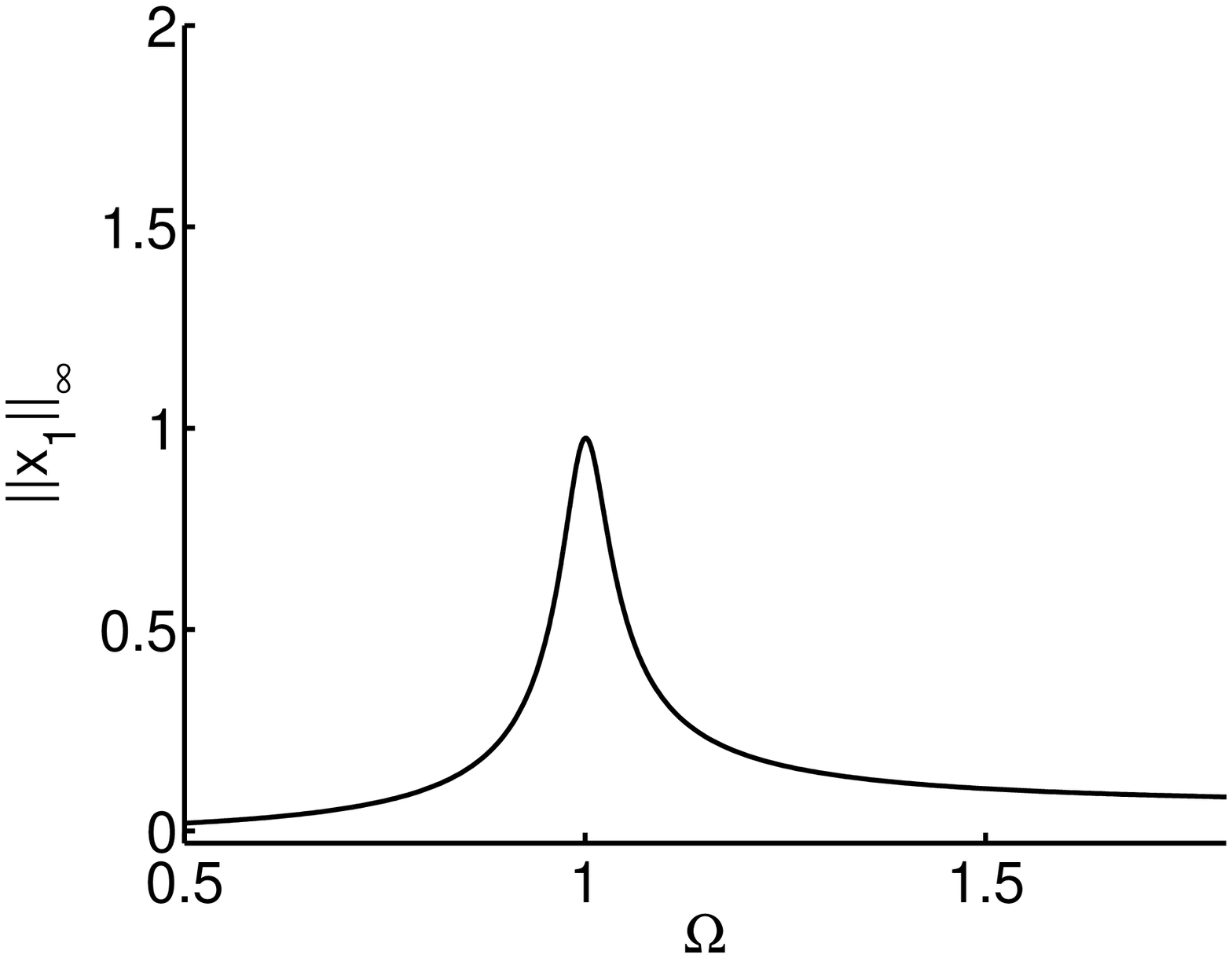}
\label{fig:isola_1}}

\subfloat[Subfigure 3 list of figures text][]{
\includegraphics[width=0.45\textwidth]{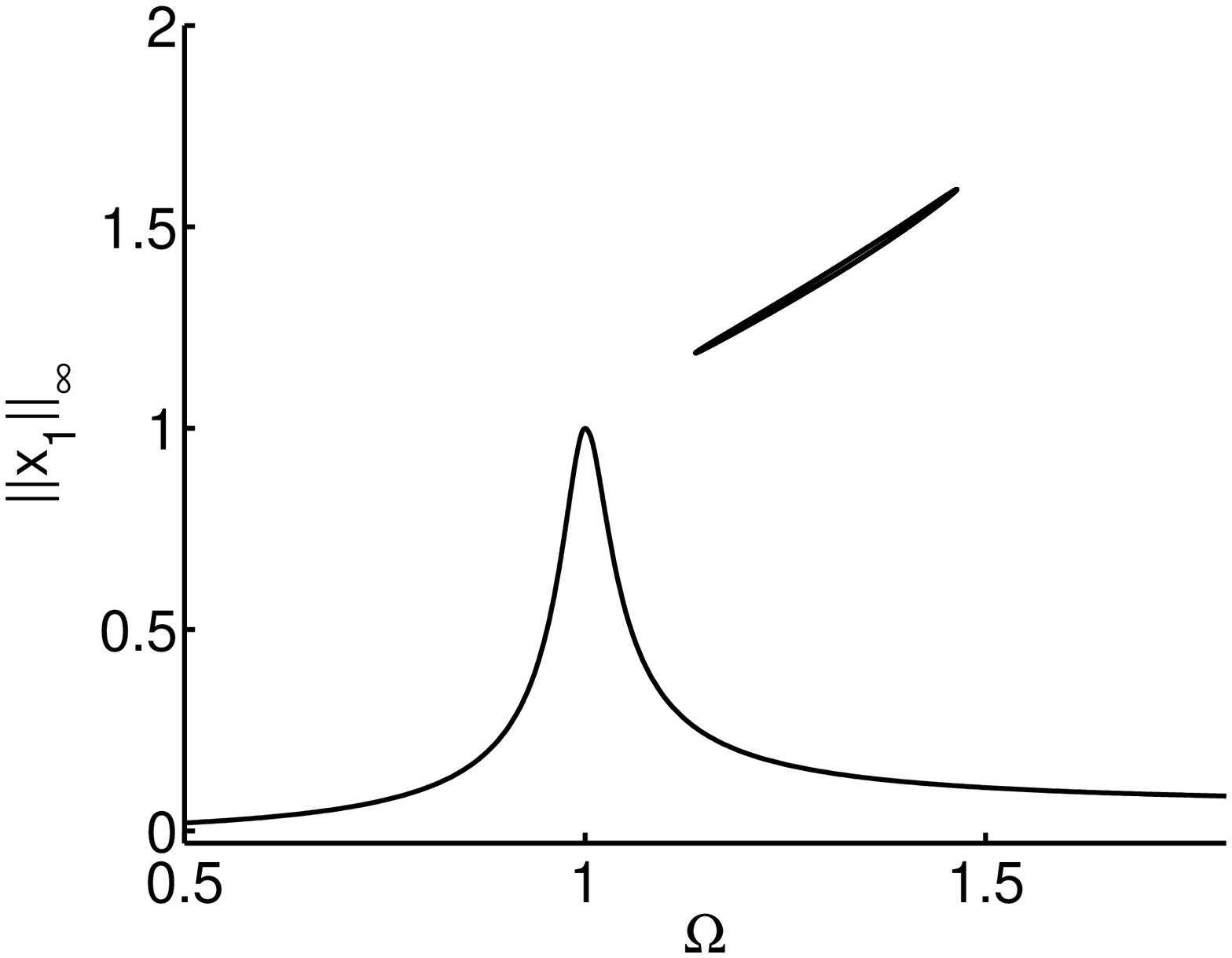}
\label{fig:isola_2}}
\hfill
\subfloat[Subfigure 3 list of figures text][]{
\includegraphics[width=0.45\textwidth]{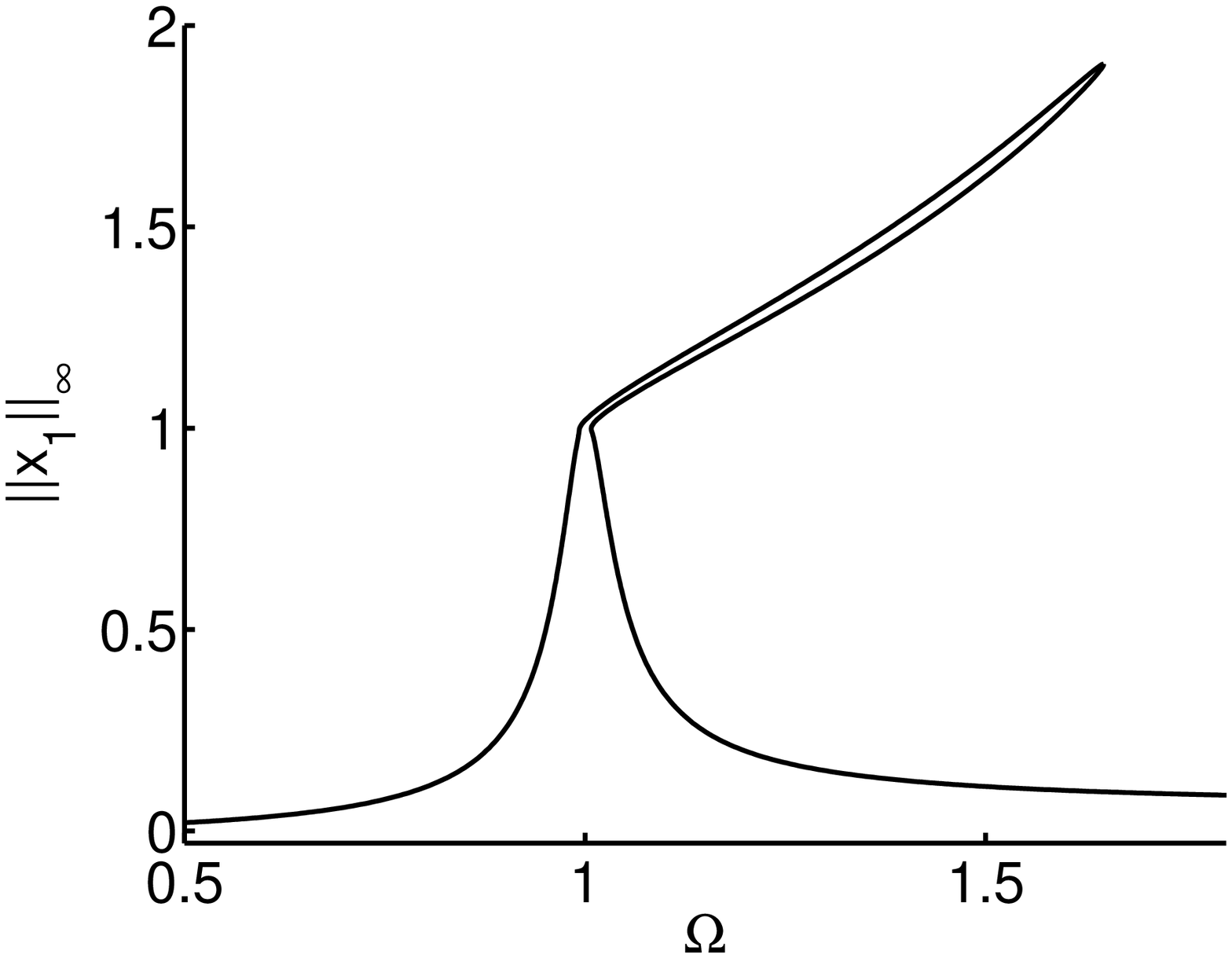}
\label{fig:isola_3}}

\caption{The solid curve in Panel \protect\subref{fig:lp_isola_curve} is a projection of the loci of fold points, while the three horizontal lines correspond to the sections for $I_l\in\{0.0463,0.0475,0.0487\}$, which are shown in Panels \protect\subref{fig:isola_1}--\protect\subref{fig:isola_3}, respectively. Panel \protect\subref{fig:isola_1} shows the frequency response of $I_l=0.0463$; panel \protect\subref{fig:isola_2} shows the isola for the frequency response of $I_l=0.0475$; panel \protect\subref{fig:isola_3} shows the frequency response of $I_l=0.0487$. The remaining parameter set is $(\xi,\beta)=(0.03,0.885)$.}
\label{fig:isola}
\end{figure}

In addition to the difference in the number of fold points between the case of continuous vs.\ discontinuous forcing amplitude the latter exhibits an isola of periodic solutions for all sections in a small range of the rescaled forcing amplitude $I_l=\frac{d}{m}\frac{A}{\Delta_l}$; where $\Delta_l$ denotes the displacement of the end-point of the beam when the beam touches the mechanical stops. This isola is not connected to the main solution branch. The structure of the branching process of the isola is shown in Figure \ref{fig:isola}; where Figure \ref{fig:isola}\protect\subref{fig:lp_isola_curve} is a projection of the fold point curve that connects the two right-most fold points in the sections showing the nonlinear resonance tongue, e.g.,\ Figure \ref{fig:globfig03}\protect\subref{fig:alpha5nu1_freqresp}. Note that the fold point curve that connects the two left-most fold points is not present because the forcing amplitude is to small. The horizontal lines are the sections which are shown in Figures \ref{fig:isola}\protect\subref{fig:isola_1}--\ref{fig:isola}\protect\subref{fig:isola_3} and combined they illustrate structure of the isola. In Figure \ref{fig:isola}\protect\subref{fig:isola_1} the section with $I_l=0.0463$ is shown and it is readily observed that the frequency response exhibits only the main solution branch for the system \eqref{eq:BD_system} with no impacting solutions, i.e.,\ a linear resonance peak as a result of weak forcing; in Figure \ref{fig:isola}\protect\subref{fig:isola_2} the section with $I_l=0.0475$ is shown and an isola of periodic solutions appears; increasing the forcing amplitude to  $I_l=0.0487$ and the isola has reconnected to the main solution branch. In the present case the isola is caused by the discontinuous change in forcing amplitude. In Figure \ref{fig:waterfall}\protect{\subref{fig:waterfall_surf}} a global two-parameter bifurcation diagram is shown, this will be explained in detail in Section \ref{sec:exp_comp}.\par

\section{Comparison with experiments}\label{sec:exp_comp}

The experimental data is adapted from \cite{Bureau2013}, where experimental continuation is performed for an impacting beam experiment. A cantilever beam of length $0.161\, {\rm m}$ with end mass $0.2116\, {\rm kg}$ is periodically excited by a shaker and allowed to hit stops at $0.71\, {\rm m}$ with a gap between cantilever and stops of $10^{-3}\, {\rm m}$ (cf.\ Figure \ref{fig:system_sketch}). The data set was obtained via frequency sweeps with constant gains of the shaker. As a result the amplitudes of the excitation are not constant along the constant-gain sweeps conducted in the experiment. We transformed this data set for comparison with the numerical findings of the model \eqref{eq:BD_system}. The non-constant forcing amplitudes during sweeps prohibit conclusive evidence of the presence of isolas in the experiment, as observed in the model; the necessary information is simply not contained in this data set. The coupling of the shaker to the mechanical structure resulted in some cases in forced excitations that are not perfect harmonics. This could be the reason why, in the experimental data the cusp point also comes earlier than in the model. We now consider system \eqref{eq:BD_system} with $\xi, \beta, \alpha, \nu$ and $p$ are fixed, while $I_l$ and $\Omega$ are free parameters. The values used are $\xi = 0.03, \beta=0.885, \alpha=5.9,\nu=1\text{ and }p=10^2$. The experimentally fitted value of the change in forcing amplitude was found to be $\nu=1$, even though the theoretical predictions indicate that $\nu$ should be of higher order and, therefore, negligible.
\begin{figure}[h]
\centering

\subfloat[Subfigure 1 list of figures text][]{
\includegraphics[width=0.9\textwidth]{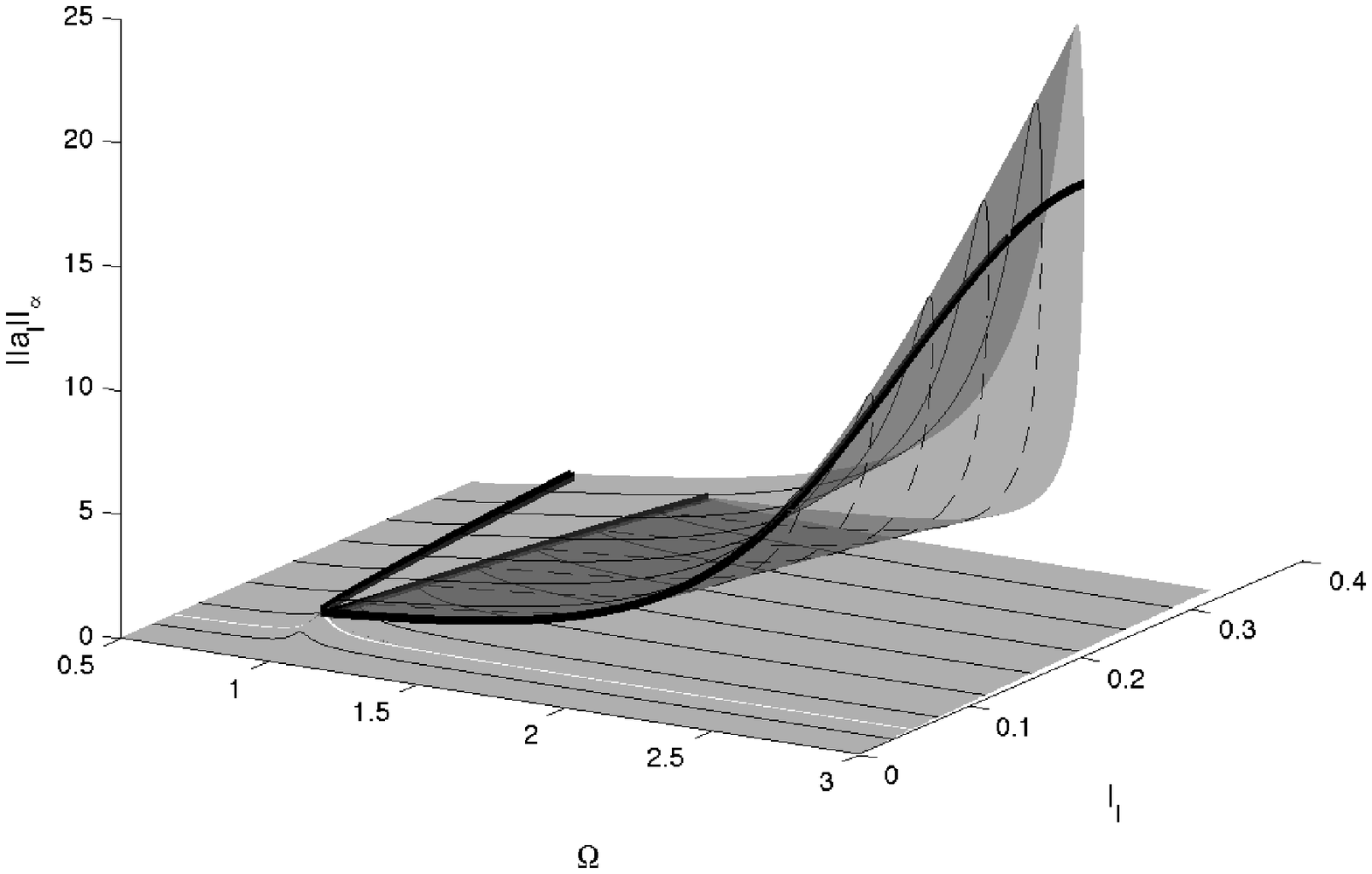}%
\label{fig:waterfall_surf}}

\subfloat[Subfigure 2 list of figures text][]{
\includegraphics[width=0.45\textwidth]{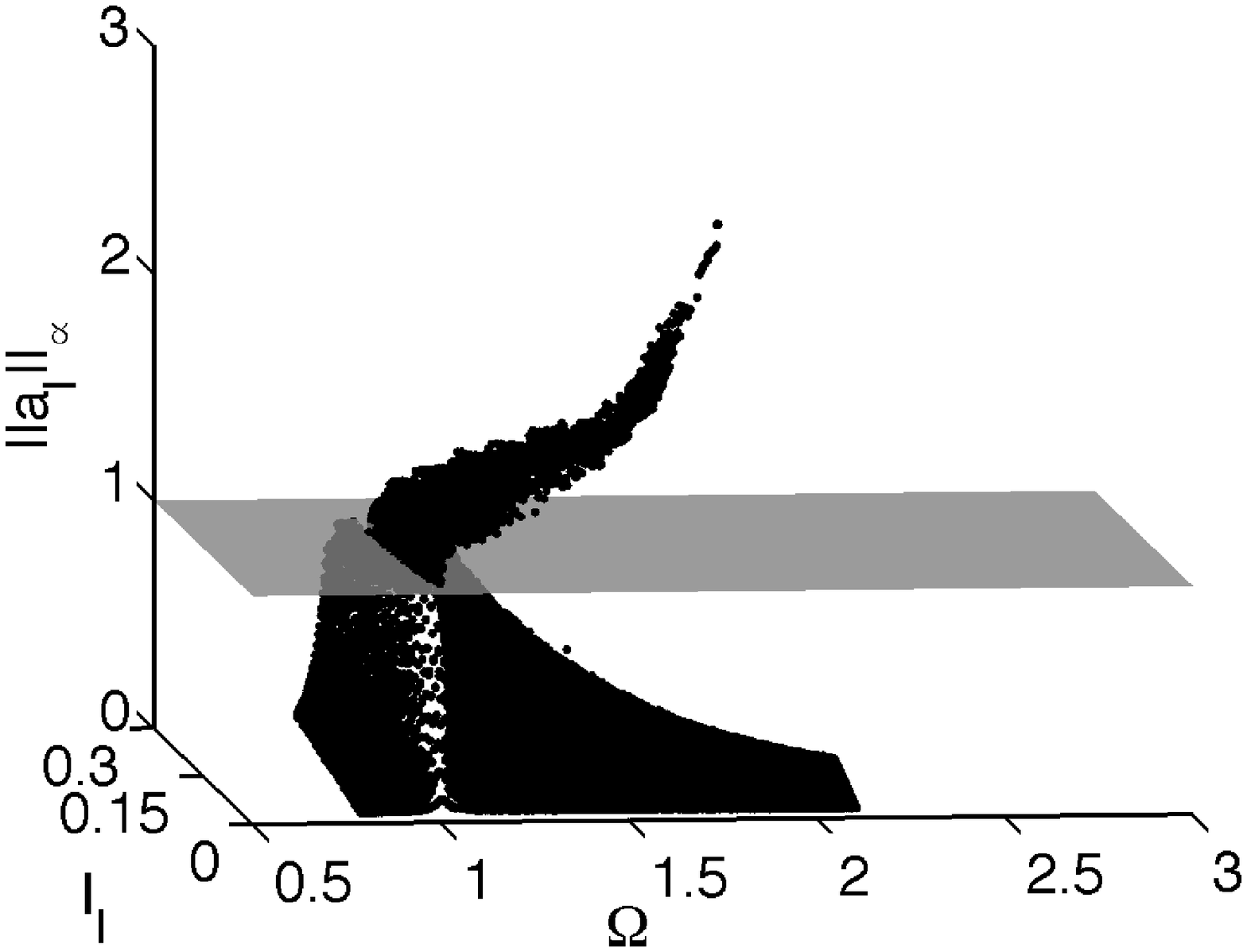}
\label{fig:impact_surf_exp}}\hfill
\subfloat[Subfigure 3 list of figures text][]{
\includegraphics[width=0.45\textwidth]{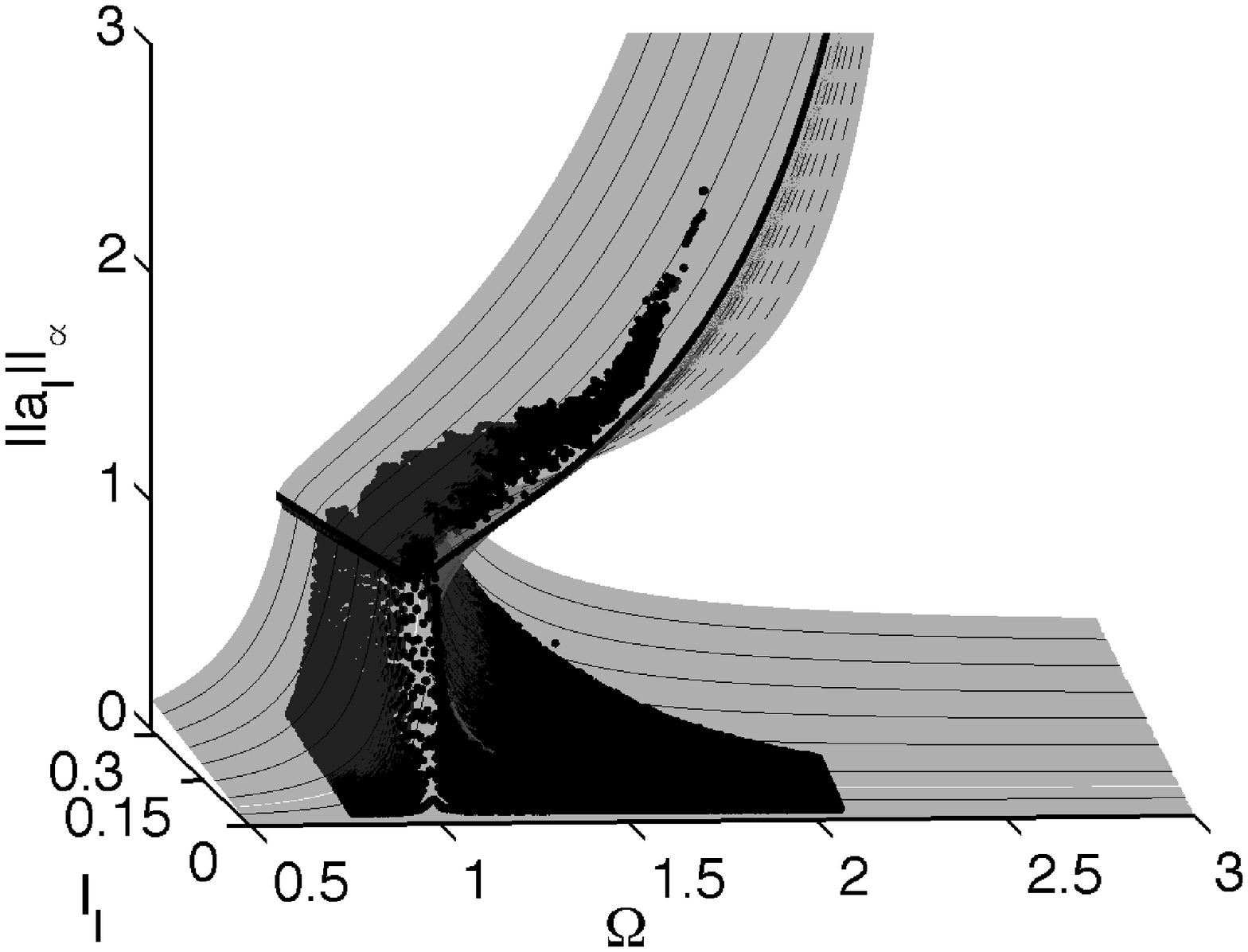}
\label{fig:waterfall_surf_exp}}

\caption{Panel \protect\subref{fig:waterfall_surf} shows the global bifurcation structure of the 1:1 resonance tongue in two parameters, namely, driving force amplitude $I_l$ and driving frequency $\Omega$. The thin curves will be used for experimental comparison. The solid black curves represent the loci of fold points. Panel \protect\subref{fig:impact_surf_exp} shows the experimental data plotted against the surface indicating the mechanical stops. Panel \protect\subref{fig:waterfall_surf_exp} shows the smooth model \eqref{eq:BD_system} with the experimental data adapted from \cite{Bureau2013}.}
\label{fig:waterfall}
\end{figure}
\clearpage

We are interested in how the loci of fold points depend on driving frequency $\Omega$ and the rescaled forcing amplitude $I_l$. From now on we also consider a rescaling of the displacement of the end point of the beam by $\Delta_l$; this is denoted by $a_l$ and $||a_l|_{\infty}$ is the maximal displacement amplitude of the periodic solutions. The scaling is introduced because the laser in the experiment measures the displacement at a certain distance below the lumped mass and for experimental comparison $a_l$ is adjusted accordingly using the beam configuration in Figure \ref{fig:shapes} and assuming small displacements; this transformation has no effect on qualitative measures. The global two-parameter bifurcation diagram is shown in Figure \ref{fig:waterfall}\protect{\subref{fig:waterfall_surf}}; here we see how the resonance tongue is bending over due to the mechanical stops and the nonlinearity they induce. In the diagram there are two solid black curves, separated via the plane section given by \mbox{$\lbrace(\Omega,I_l,||a_l||_\infty) : \Omega=1\rbrace$}, and these are the loci of fold points and both of the curves loop back over themselves in a cusp in the projection onto the ($\Omega; I_l$)-parameter plane. The range over which the isola exists is visualized by cutting the manifold in two pieces leaving out the sections with the isola. This isola is very thin and can not be observed in the experimental data. The thin solid/dashed curves are used for experimental comparison later.\par
\begin{figure}[h]
\centering
\subfloat[Subfigure 1 list of figures text][]{
\includegraphics[width=0.45\textwidth]{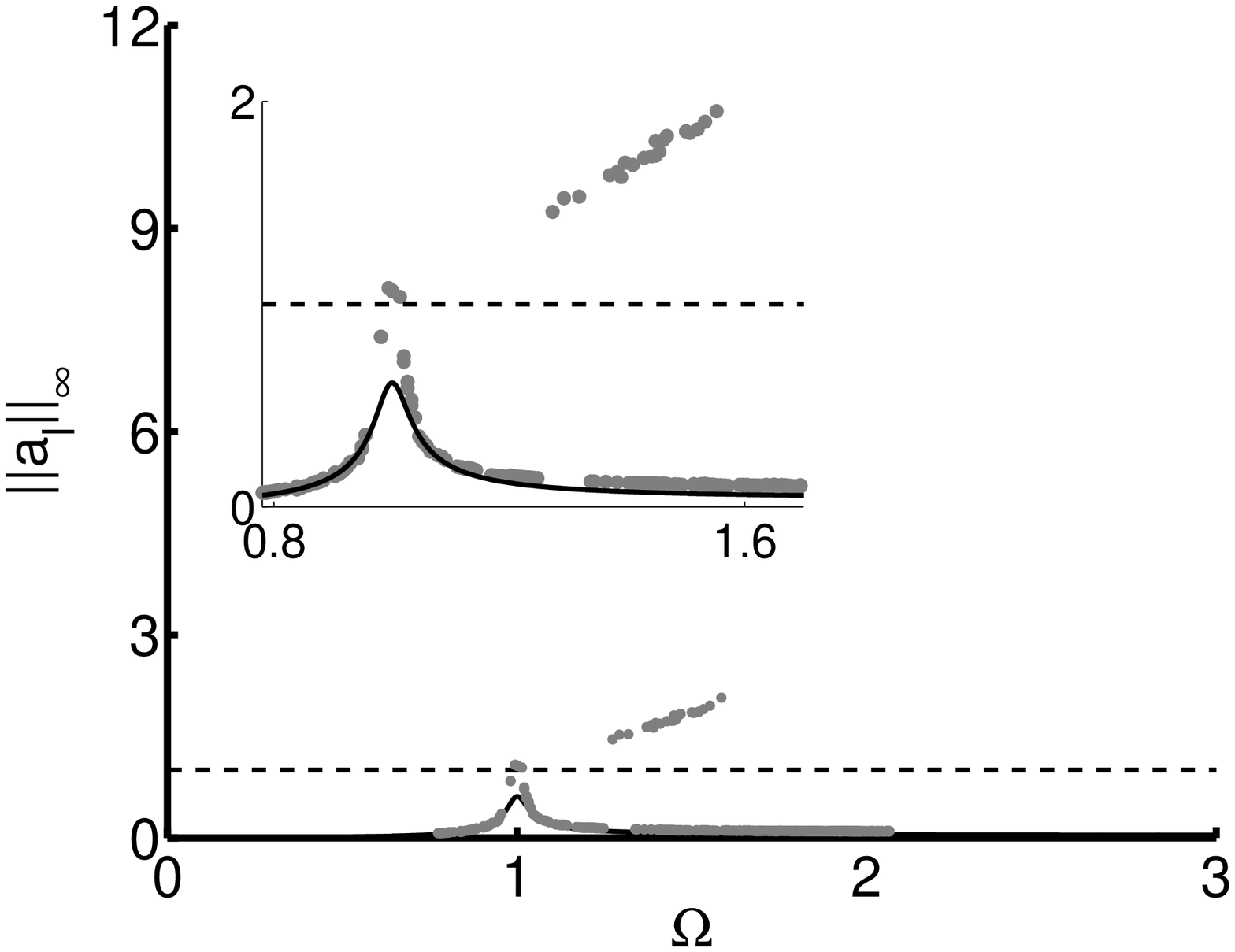}
\label{fig:exp_comp_idA_exp2}}\hfill
\subfloat[Subfigure 2 list of figures text][]{
\includegraphics[width=0.45\textwidth]{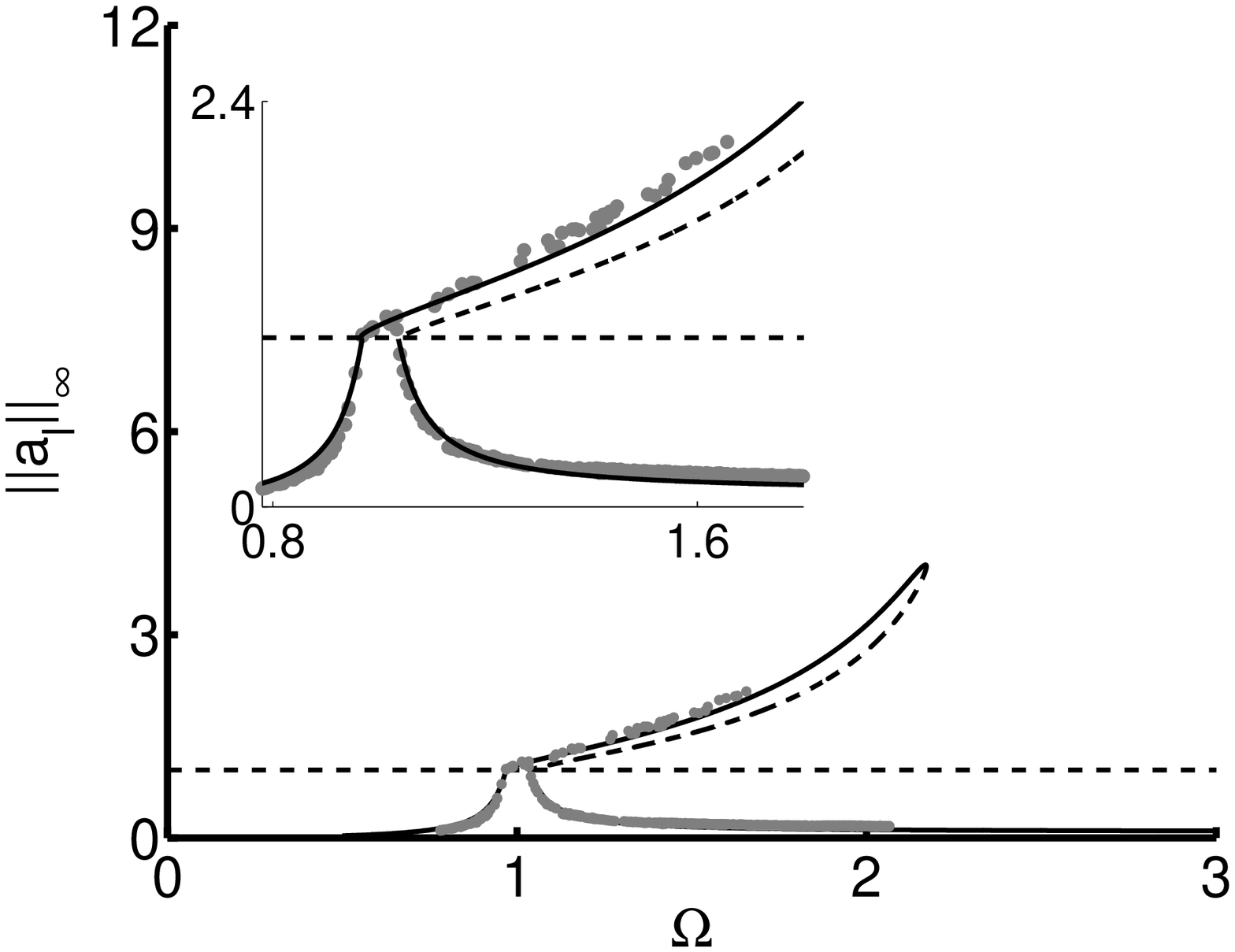}
\label{fig:exp_comp_idA_exp4}}
\vfill
\subfloat[Subfigure 3 list of figures text][]{
\includegraphics[width=0.45\textwidth]{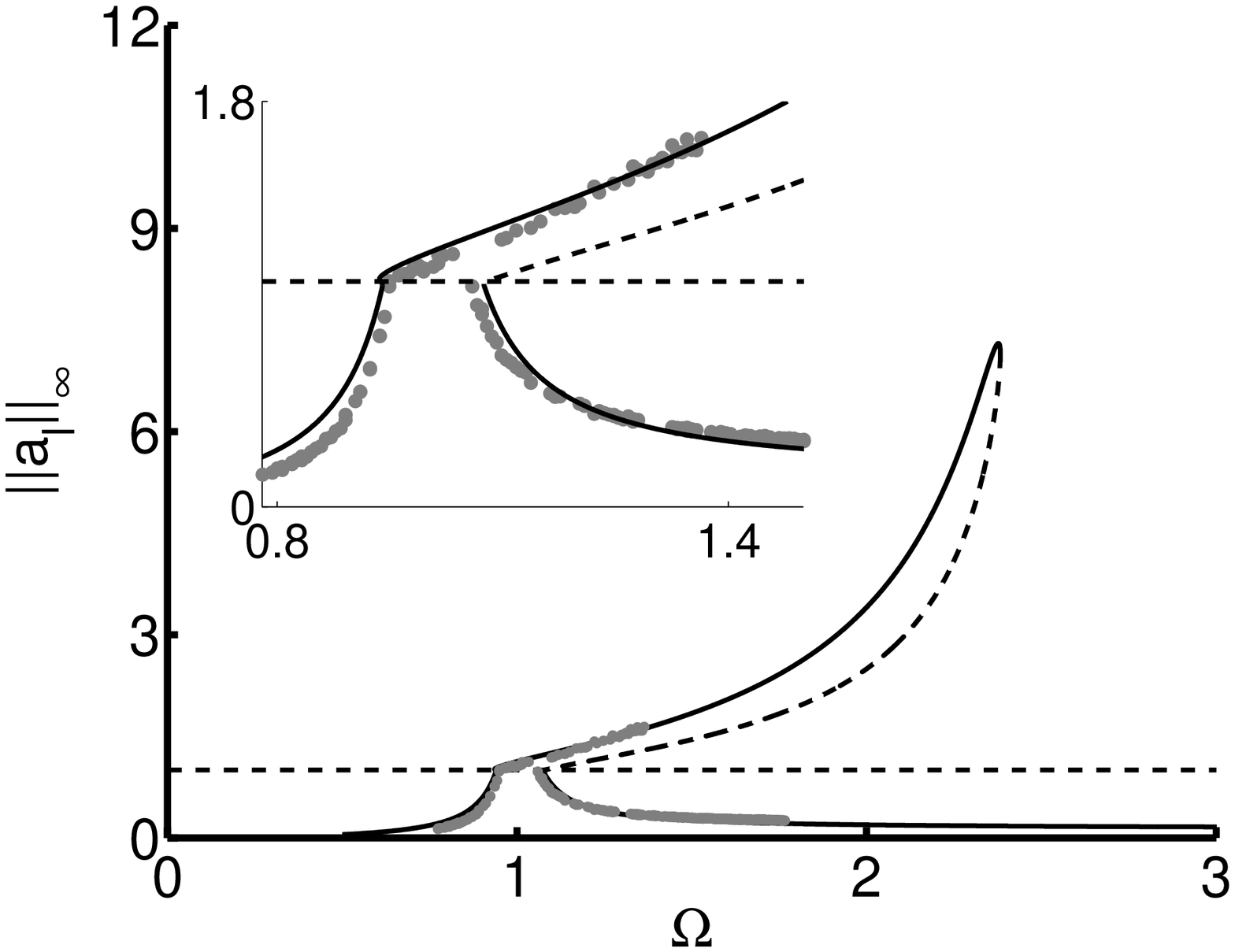}
\label{fig:exp_comp_idA_exp6}}\hfill
\subfloat[Subfigure 4 list of figures text][]{
\includegraphics[width=0.45\textwidth]{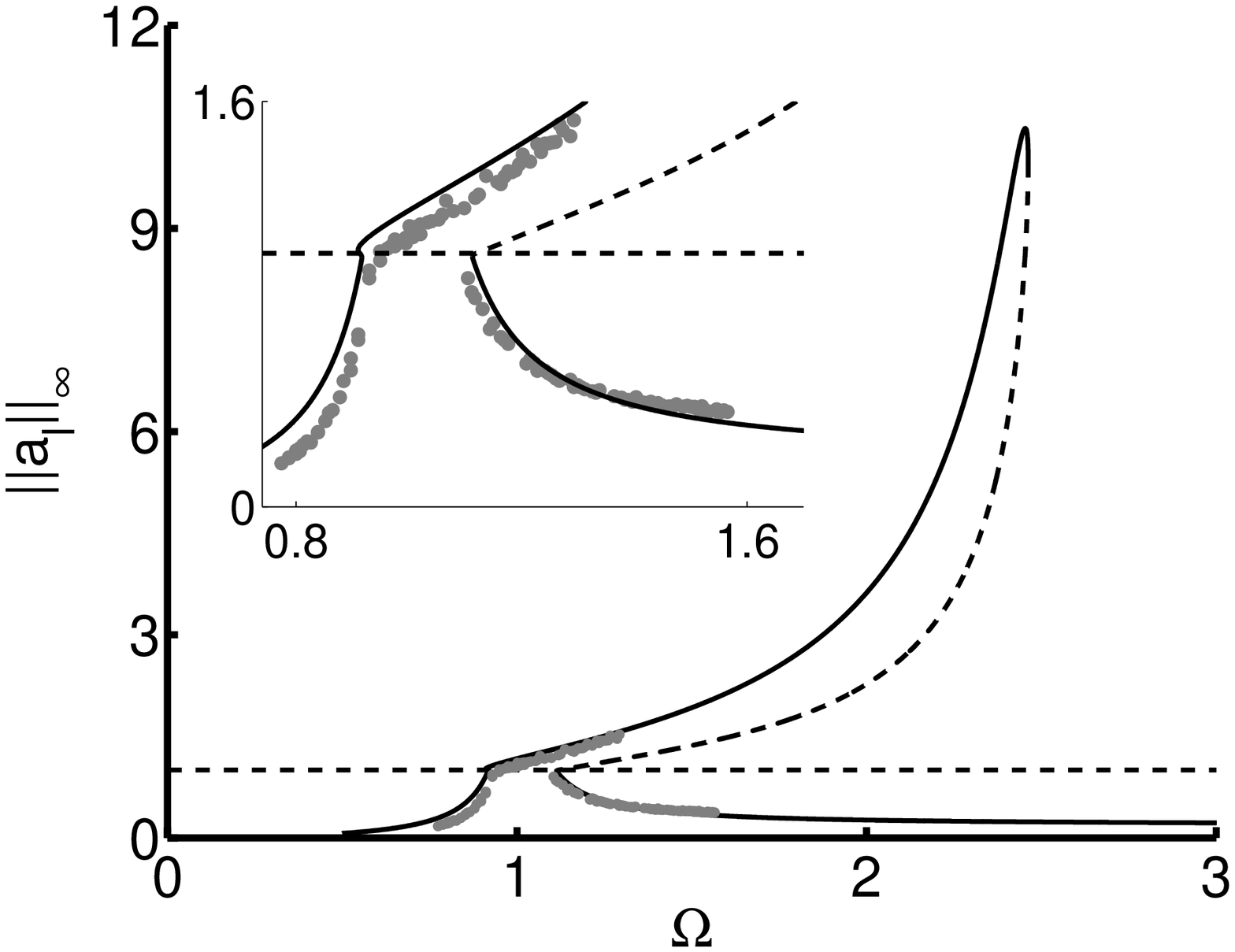}
\label{fig:exp_comp_idA_exp8}}
\vfill
\subfloat[Subfigure 4 list of figures text][]{
\includegraphics[width=0.45\textwidth]{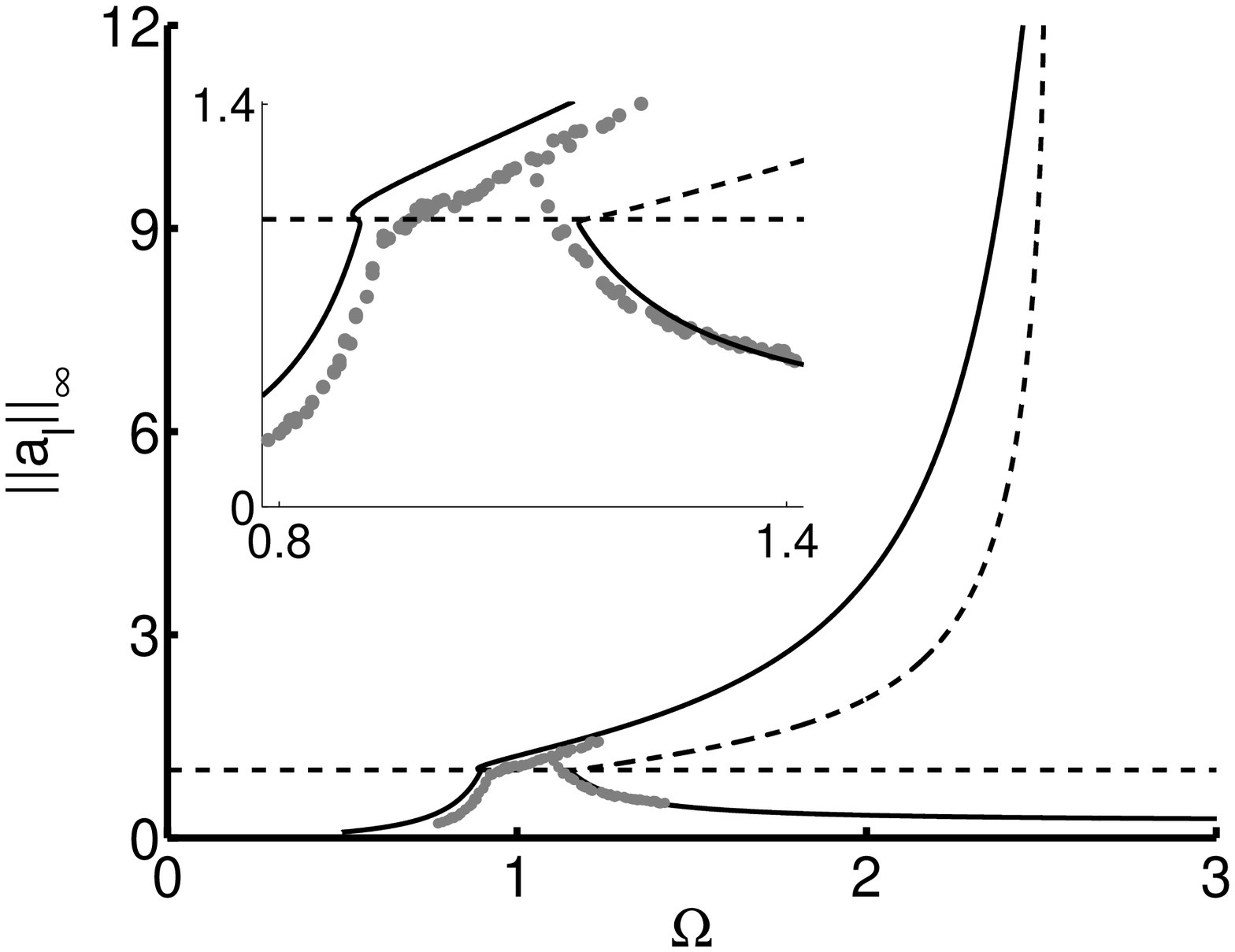}
\label{fig:exp_comp_idA_exp10}}\hfill
\subfloat[Subfigure 4 list of figures text][]{
\includegraphics[width=0.45\textwidth]{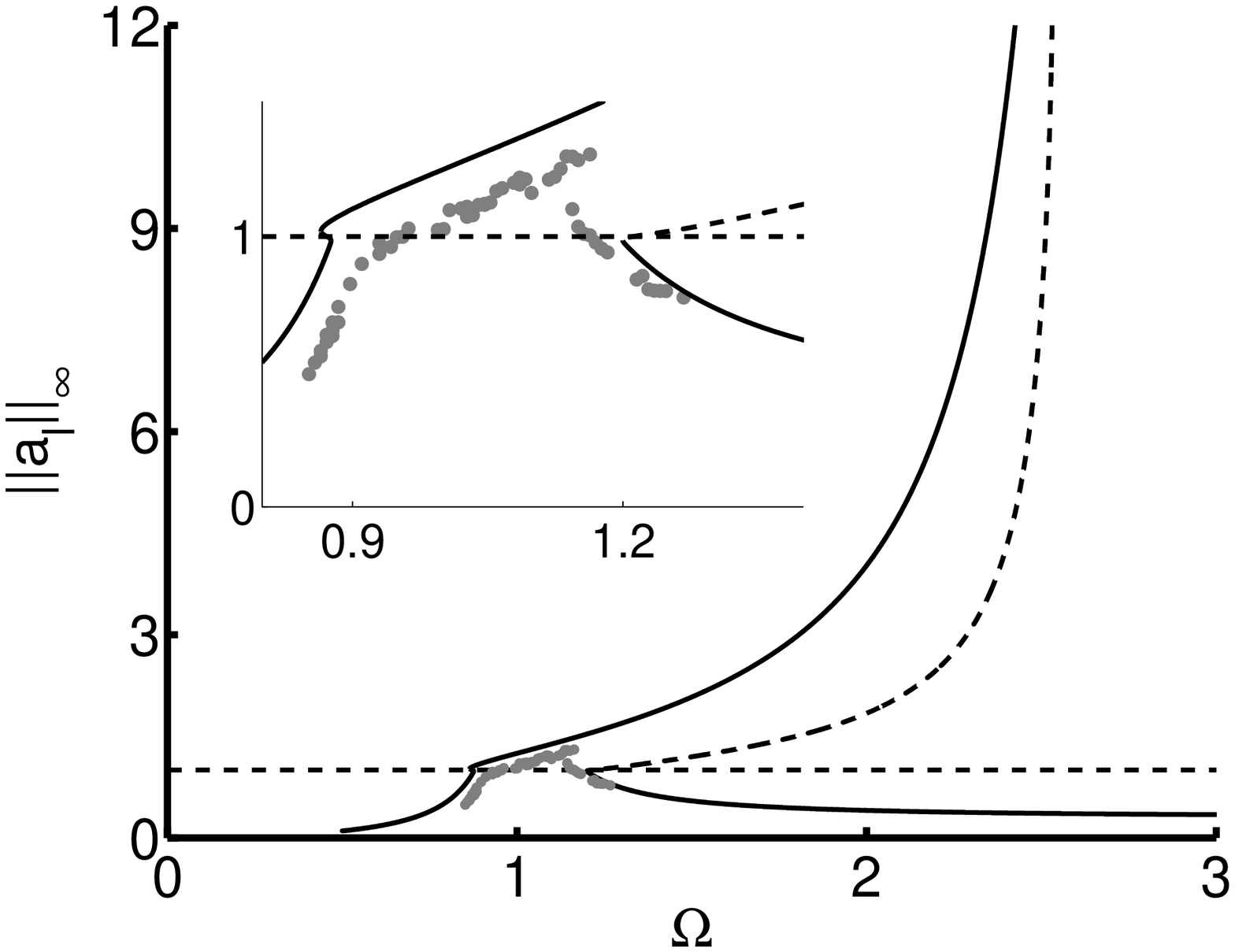}
\label{fig:exp_comp_idA_exp12}}
\caption{Panels \protect\subref{fig:exp_comp_idA_exp2}-\protect\subref{fig:exp_comp_idA_exp12} shows the comparison of the smooth model \eqref{eq:BD_system} with the experimental data set. Each plot represents one of the thin solid/dashed curve sections in Figure \ref{fig:waterfall}. The smoothing homotopy parameter is $p=100$ and $I_l$ takes the values 0.03, 0.07, 0.11, 0.16, 0.20 and 0.24, respectively. The solid/dashed curve represents the continuation results of \eqref{eq:BD_system} and the grey bullets are the scattered point set from the experimental data. The enlargements are for proper visual comparison in the regions with experimental data.}
\label{fig:exp_comp_slices}
\end{figure}

 
The experimental data set is shown in Figure \ref{fig:waterfall}\protect{\subref{fig:impact_surf_exp}} with the mechanical stop indicated by the transparent surface, and we can clearly see the match of the resonance tongue bending over due to the mechanical stops. Figure \ref{fig:waterfall}\protect{\subref{fig:waterfall_surf_exp}}  shows an enlargement of the continuation results from the smooth system \eqref{eq:BD_system} in $(\Omega,I_l,||a_l||_{\infty})$ together with the experimental data. It is difficult to asses the quantitative agreement between the model and experiments because the data is scattered and the density of measurements varies nonuniformly. Therefore, all experimental comparison is done via sections of constant forcing shown in Figure \ref{fig:exp_comp_slices}.\par 

Figure \ref{fig:exp_comp_slices}\protect\subref{fig:exp_comp_idA_exp2} shows the section for $I_l\approx 0.03$. The model predicts low-amplitude nonimpacting solutions for this value of (scaled) forcing amplitude, but the experimental data is showing impacting solutions. This indicates that the cusp point or grazing point is not predicted perfectly. From the figure it looks as if the experimental data indicates an isola as well. As mentioned earlier, it is not possible to conclude if this is indeed the case because the experimental sweeps were conducted with constant gain of the shaker. Furthermore, for small-amplitude forcing the noise levels might not be negligible.\par

Figure \ref{fig:exp_comp_slices}\protect\subref{fig:exp_comp_idA_exp4} shows the section for $I_l\approx 0.07$. For this value of the forcing amplitude the quantitative agreement of both upper and lower branches is very good. It is worth noting that, for orbits close to the mechanical stops, there is, in some of the figures, a small indication of clustering points. They could be ascribed to the fact that experimental measurements naturally show some fluctuations in regions where noise levels are comparable to the distance from the mechanical stops, and grazing orbits are particularly vulnerable to this effect.\par

Figures \ref{fig:exp_comp_slices}\protect\subref{fig:exp_comp_idA_exp6} and \ref{fig:exp_comp_slices}\protect\subref{fig:exp_comp_idA_exp8} show the sections for $I_l\approx 0.11$ and $I_l\approx 0.16$, respectively. Again, the quantitative agreement of both upper and lower branches is very good. From the enlargement it is seen that the width of the 'linear' resonance in the model is not in exact agreement, with the experimental data, and the upper branch is slightly overestimated. The reasons for these small discrepancies could be the choice of damping model or damping parameters, which are phenomenological.\par

In Figures \ref{fig:exp_comp_slices}\protect\subref{fig:exp_comp_idA_exp10} and \ref{fig:exp_comp_slices}\protect\subref{fig:exp_comp_idA_exp12}  sections are shown with $I_l\approx 0.20$ and $I_l\approx 0.24$, respectively. For these values of the forcing amplitude, the quantitative agreement is less good; in particular, the effect of overestimating the upper branches is slightly more pronounced, but it should be noted that the absolute errors are $<1\, {\rm mm}$. Furthermore, from the enlargements it can also be observed that the fold points of the lower branches are not sharply defined for these large forcing amplitudes.\par
\clearpage

The correspondence between the experimental data and the mathematical model also depends on the ambiguous choice of the ratio $\frac{m}{d}$ for the rescaling $I_l=\frac{d}{m}\frac{A}{\Delta_l}$. In these calculations $m/d=2/3$ was used and the measured distance $\delta$ in Figure \ref{fig:system_sketch} is corrected to fit by setting $\delta = 0.85\cdot10^{-3}\, {\rm m}$. Furthermore, it is to be expected that the model validity is restricted to small-amplitude displacement, because for larger displacement amplitudes the beam should be treated as a nonlinear elastica.\par 
It would have been desirable to investigate whether or not the two limit points before the first linear resonance exist in the experiment or not. Unfortunately, this is not feasible, because these two points lie at a distance below the noise level of the experiment.\par 

\section{Discussions and conclusions}
In this paper a single-degree-of-freedom model of a cantilever beam with symmetric mechanical stops and a mass attached was derived. Model parameters were estimated from first principles and compared to experimentally fitted values. These comparisons are given in Table \ref{tab:Parameter_estimation} and it should be remarked that the nonlinear stiffness ratio $\alpha$ is quite well approximated. The model is particularly suited for small-amplitude solutions, for which it predicts the experimental data very well.\par
In contrast to the more traditional setting where mechanical systems with impacts are realized by nonsmooth models, the present paper dealt with the setting in which the impact model is described with a smooth model. To obtain this model a smooth switching (nonlinear homotopy) was applied. The smoothing procedure was investigated by means of numerical bifurcation analysis. It was investigated when and how the bifurcation structure was affected by the smoothing. The investigation showed that, as the smoothing function approaches the Heaviside function, the bifurcation diagram becomes practically independent of the smoothing homotopy parameter. This investigation was performed for moderate impacts, hard impacts and discontinuous forcing amplitudes. The mechanism for the smoothing-induced bifurcations was also identified, and it was shown that the suitable smoothing homotopy parameter depends on the stiffness ratio, i.e.,\ $p=\Theta(\alpha)$; specifically, increasing $\alpha$ implies increasing $p$. With this in mind, it may be noted that the smoothing process can cause numerical difficulties as $\alpha$ becomes very large, because then $p$ must be very large. In such cases it may be beneficial to perform multisegment numerical continuation of the nonsmooth system.\par
The behavior of the smooth model was compared with the experimental data \cite{Bureau2013}. The qualitative structure was verified within the range of the data set. The quantitative comparison was good overall, with a tendency to overestimate the upper branch as the forcing amplitudes becomes large. The results in the paper demonstrate that a smooth single-degree-of-freedom model describes the experimental behavior very well and that high-dimensional FE models are not always necessary to obtain predictions of good quantitative quality.\par
For future work it is interesting to test the model in a larger parameter regime to explore its limits. Specifically it would be interesting to investigate other parameter regimes of excitation frequency and amplitude. Since a posteriori fitting was kept to a minimum, it could also be interesting to find an optimal fit in a reasonable parameter-neighborhood, now that model parameters have been determined with mechanical interpretations and proper magnitudes.

\section*{Acknowledgements}

ME thanks the Department of Mathematics at \emph{The University of Auckland} for its kind hospitality and the Idella Fondation for financial support. JS and JJT acknowledge funding from the Danish Research Council FTP under the project number 09-065890/FTP. The authors thank Emil Bureau, Frank Schilder and Ilmar Santos for fruitful discussions and for providing the data from \cite{Bureau2013}.

\bibliographystyle{unsrt}

\bibliography{references}

\end{document}